\documentclass[11pt]{article}

\usepackage{xcolor}
\usepackage{siunitx}
\usepackage[export]{adjustbox}
\usepackage{authblk}
\usepackage{tabularx}
\usepackage{graphicx}
\usepackage{paralist}
\usepackage{url}
\usepackage{wrapfig}
\usepackage{amsmath}
\usepackage{amssymb}
\usepackage{placeins}
\usepackage{tikz, standalone}
\usepackage{caption}
\usepackage{subcaption}

\usetikzlibrary{shapes,positioning,intersections,quotes, calc, angles}
\usepackage{color}

\newcommand{\sfrac}[2]{\mathchoice
  {\kern0em\raise.5ex\hbox{\the\scriptfont0 #1}\kern-.15em/
   \kern-.15em\lower.25ex\hbox{\the\scriptfont0 #2}}
  {\kern0em\raise.5ex\hbox{\the\scriptfont0 #1}\kern-.15em/
   \kern-.15em\lower.25ex\hbox{\the\scriptfont0 #2}}
  {\kern0em\raise.5ex\hbox{\the\scriptscriptfont0 #1}\kern-.2em/
   \kern-.15em\lower.25ex\hbox{\the\scriptscriptfont0 #2}}
  {#1\!/#2}}

\newcommand{\half}{{\sfrac{1}{2}}}
\newcommand{\nph}{{n+\sfrac{1}{2}}}
\newcommand{\nmh}{{n-\sfrac{1}{2}}}

\DeclareMathSymbol{\shortminus}{\mathbin}{AMSa}{"39}

\def\half   {\frac{1}{2}}

\def\nph    {n+\sfrac{1}{2}}
\def\iph    {i+\sfrac{1}{2}}
\def\jph    {j+\sfrac{1}{2}}
\def\kph    {k+\sfrac{1}{2}}
\def\imh    {i-\sfrac{1}{2}}
\def\jmh    {j-\sfrac{1}{2}}
\def\kmh    {k-\sfrac{1}{2}}

\def\uvec   {\mathbf{u}}
\def\umachi   {u^{MAC}_{\iph,j,k}}
\def\umaclo   {u^{MAC}_{\imh,j,k}}
\def\vmachi   {v^{MAC}_{i,\jph,k}}
\def\vmaclo   {v^{MAC}_{i,\jmh,k}}
\def\wmachi   {w^{MAC}_{i,j,\kph}}
\def\wmaclo   {w^{MAC}_{i,j,\kmh}}

\newcommand{\HeatFlux}{\boldsymbol{\mathcal{Q}}}
\newcommand{\SpeciesFlux}{\boldsymbol{\mathcal{F}}}

\newcommand{\StressTensor}{\boldsymbol{\tau}}

\providecommand{\keywords}[1]
{
  \small	
  \textbf{\textit{Keywords---}} #1
}

\begin{document}

\title{A Weighted State Redistribution Algorithm \\
for Embedded Boundary Grids}

\author[a]{A. Giuliani}
\author[b]{A.S. Almgren}
\author[b]{J.B. Bell}
\author[a,c]{M.J. Berger}
\author[d]{\authorcr M.T. Henry de Frahan}
\author[e,f]{D. Rangarajan}

\affil[a]{Courant Institute, NYU }
\affil[b]{Lawrence Berkeley National Laboratory}
\affil[c]{Flatiron Institute, NYC}
\affil[d]{National Renewable Energy Laboratory}
\affil[e]{National Energy Technology Laboratory, Morgantown, WV  26507, USA}
\affil[f]{NETL Support Contractor, Morgantown, WV  26507, USA}

\date{}
\maketitle

\begin{abstract}
State redistribution is an algorithm that stabilizes cut cells for embedded boundary grid methods. This work 
extends the earlier algorithm in several important ways. 
First, state redistribution is extended to three spatial dimensions.  Second, we discuss several algorithmic changes and improvements motivated by the more complicated cut cell geometries that can occur in higher dimensions. In particular, we introduce a weighted version with less dissipation in an easily generalizable framework.  Third, we demonstrate that state redistribution can also stabilize a solution update that includes both advective and  diffusive contributions. 
The stabilization algorithm is shown to be effective for 
incompressible as well as compressible reacting flows.  
Finally, we discuss the implementation of the algorithm for several exascale-ready 
simulation codes based on AMReX,  demonstrating  ease of use in combination with
domain decomposition, hybrid parallelism and complex physics. 
\end{abstract}

 \section{Introduction}\label{intro}
 Embedded boundary grids are useful for solving partial differential equations on complex engineering domains since mesh generation is robust and automatic.  These grids are composed of regular Cartesian cells away from the boundary, covered cells that do not participate in the solution procedure, and irregular cut cells that intersect the embedded boundary.  It is well known that the cut cells can be arbitrarily small, thus care must be taken when using explicit time stepping methods on these grids.  Many different approaches that address the small cell problem have been proposed over the years. Most intuitive with a lot of recent progress is cell merging \cite{Saye17}. The interesting work in \cite{gokhaleNikosKlein:2018} is dimensionally split, but not that accurate at the cut cells.  $H$-box methods \cite{mjb-hel-rjl:hbox} have nice theoretical properties but seem difficult to implement in three dimensions.

 One technique that has been successfully used in three dimensions on complex geometries is called flux redistribution (FRD) \cite{chern1987conservative, COLELLA2006347}.  The idea is that each cell takes its maximum stable time step, then excess flux is redistributed to neighboring cells to maintain conservation.  This approach is attractive because it is simple to implement as a postprocessing step. However, the downside is that there is a loss of accuracy at the cut cells and the scheme is not linearity preserving.  
 Recently, a new approach called state redistribution (SRD) was proposed for finite volume methods on two-dimensional cut cell grids \cite{BG}.  This is a minimally invasive stabilization technique that is linearity preserving, conservative, and straightforward to implement in two dimensions for hyperbolic conservation laws.  Inspired by flux redistribution, SRD postprocesses the numerical solution by accurately redistributing the solution states in a way that maintains conservation. 
 
  This work builds on \cite{BG} in several important ways. The work in \cite{BG} was in two space dimensions. 
  As we show here, the extension to 3D  was straightforward.  We propose a weighted version of SRD with reduced dissipation, that smoothly activates depending on the volume fraction of a cut cell.
  The framework of this weighted algorithm will help to generalize to other types of weightings, e.g., to preserve monotonicity.
   We discuss some of the design choices made for easier implementation in a 3D parallel production code. Previous work on SRD focused on the Euler equations of gas dynamics. We show here that SRD is simple enough to apply to many types of equations. The extensions shown later  include compressible reacting flow and low Mach number multiphase flow.
  Finally, the algorithm has been incorporated into AMReX-Hydro, a set of modules based on AMReX \cite{amrex:ijhpca}. This is a software framework that supports the development of block-structured adaptive mesh refinement (AMR) algorithms for solving systems of partial differential equations on simple or complex geometries, using machines from laptops to exascale architectures.
  While all these examples use finite volume schemes on the embedded boundary mesh, we mention that SRD has also been extended to discontinuous Galerkin schemes in two dimensions \cite{giuliani2021two}, and research on higher order finite volume SRD schemes is in progress.
  
  The rest of this paper is organized as follows. In section \ref{sec:prelim} we set the stage for the geometry representation, and outline the three dimensional version of the original algorithm  in \cite{BG}. Section \ref{sec:alphabeta} discusses a generalization of SRD that smoothly activates as the volume fraction decreases below a threshold.  It is presented in a general framework which reveals the constraint that the weights must satisfy for conservation.   Section \ref{sec:par} presents implementation details in AMReX. Section \ref{sec:algext} discusses some particular choices that go into the extensions to other sets of equations.  Computational examples are in section \ref{sec:compResults}. Conclusions and directions for future research are  in section \ref{sec:conclusions}.

\section{Preliminaries}\label{sec:prelim}
We briefly describe the the embedded boundary representation and geometry generation to give the reader an idea of the context in which State Redistribution is applied. We then present the extension to 3D of the original SRD algorithm, both the pre- and post- processing. This will provide a basis for comparison when the new weighted version is presented in section \ref{sec:alphabeta}.

\subsection{Embedded Boundary Data Structures}
When an embedded boundary (EB) is present in the domain, there are three types of cells: regular, cut and covered.  AMReX provides data structures for
accessing EB information, which is precomputed and stored in a distributed
database at the beginning of the calculation.  It is available for AMR
meshes at each level and for coarsened meshes in multigrid solvers.  The
information includes cell type, cell centroid, volume fraction,  face area
fractions and face centroids.  For cut cells, the information also
includes the centroid, normal and area of the EB face.  There is at most one EB face per cell; this is enforced at the time the geometry is generated.
Additionally, there is connectivity information between neighboring
cells.

\subsection{Geometry Generation}
AMReX uses an implicit function approach for generating the
geometry information.  The implicit function
 describes the surface of the embedded object.  It returns a
positive value, a negative value or zero, for a given position inside
the body, inside the fluid, or on the boundary, respectively.
Implicit functions for various simple shapes such as boxes, cylinders,
spheres, etc., as well as a spline based approach, are provided.
Furthermore, basic operations in constructive solid geometry (CSG) such as
union, intersection and difference are used to combine objects
together.  Geometry transformations (e.g., rotation and translation)
can also be applied to these objects. In addition to an implicit function, 
an application code can also use its own approach to
generate the geometric information and store it in AMReX's EB database.
In the current AMReX mesh generator, split cells and tunnel cells are not supported. Multivalued cells and edges with multiple cuts are not allowed.

\subsection{3D SRD Algorithm}\label{sec:srdAlg}
 We describe here the original SRD algorithm extended to  three dimensions.  The algorithm comprises a mesh preprocessing step and a solution postprocessing step. The next section will introduce some improvements and generalizations.

\subsubsection{Preprocessing}

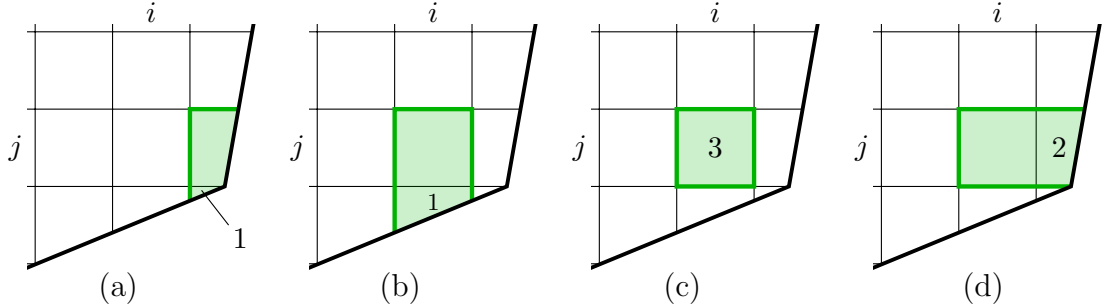
\begin{figure}
\centering
 \begin{tikzpicture}
\node (11) at (0, 0.175){  \begin{tikzpicture}

\node (a) at (0.75,2.5) {$j$};
\node (b) at (2.5,4.25) {$i$};

\clip (0.9,0.9) rectangle (4.1,4.1);
\draw (-5,-5) grid (5,5);
\draw [ultra thick, draw=green!70!black, fill=green!70!black, fill opacity=0.2] (3,1.) rectangle  (4,3);

\begin{scope}[rotate=22]
\path [fill=white] (0,.55) rectangle (5,-5);
\end{scope}
\begin{scope}[rotate=80]
\path [fill=white] (2,-3.05) rectangle (6,-6);
\end{scope}
\draw[line width=0.5mm] (0,0.58) -- (3.45,2) -- (4.45,4.2+3.45);

\node (d) at (3.65,1.325) {$1$};
\draw (3.5,1.5) -- (3.15,1.95);
  \end{tikzpicture}};
\node at (-0.25,-1.9) {\textcolor{black}{\large (a)}};

\node (11) at (3.75, 0.175){  \begin{tikzpicture}

\node (a) at (0.75,2.5) {$j$};
\node (b) at (2.5,4.25) {$i$};

\clip (0.9,0.9) rectangle (4.1,4.1);
\draw (-5,-5) grid (5,5);
\draw [ultra thick, draw=green!70!black, fill=green!70!black, fill opacity=0.2] (2,1.) rectangle  (3,3);

\begin{scope}[rotate=22]
\path [fill=white] (0,.55) rectangle (5,-5);
\end{scope}
\begin{scope}[rotate=80]
\path [fill=white] (2,-3.05) rectangle (6,-6);
\end{scope}
\draw[line width=0.5mm] (0,0.58) -- (3.45,2) -- (4.45,4.2+3.45);

\node (d) at (2.5,1.8) {\footnotesize $1$};

  \end{tikzpicture}};
\node at (3.5,-1.9) {\textcolor{black}{\large (b)}};

\node (11) at (7.5, 0.175){  \begin{tikzpicture}

\node (a) at (0.75,2.5) {$j$};
\node (b) at (2.5,4.25) {$i$};

\clip (0.9,0.9) rectangle (4.1,4.1);
\draw (-5,-5) grid (5,5);
\draw [ultra thick, draw=green!70!black, fill=green!70!black, fill opacity=0.2] (2,2.) rectangle  (3,3);

\begin{scope}[rotate=22]
\path [fill=white] (0,.55) rectangle (5,-5);
\end{scope}
\begin{scope}[rotate=80]
\path [fill=white] (2,-3.05) rectangle (6,-6);
\end{scope}
\draw[line width=0.5mm] (0,0.58) -- (3.45,2) -- (4.45,4.2+3.45);

\node (c) at (2.5,2.5) {$3$};

  \end{tikzpicture}};
\node at (7.25,-1.9) {\textcolor{black}{\large (c)}};

\node (11) at (11.25, 0.175){  \begin{tikzpicture}

\node (a) at (0.75,2.5) {$j$};
\node (b) at (2.5,4.25) {$i$};

\clip (0.9,0.9) rectangle (4.1,4.1);
\draw (-5,-5) grid (5,5);
\draw [ultra thick, draw=green!70!black, fill=green!70!black, fill opacity=0.2] (2,2.) rectangle  (4,3);

\begin{scope}[rotate=22]
\path [fill=white] (0,.55) rectangle (5,-5);
\end{scope}
\begin{scope}[rotate=80]
\path [fill=white] (2,-3.05) rectangle (6,-6);
\end{scope}
\draw[line width=0.5mm] (0,0.58) -- (3.45,2) -- (4.45,4.2+3.45);

% \node (c) at (2.5,2.5) {$3$};
\node (d) at (3.3,2.5) {$2$};

  \end{tikzpicture}};
\node at (11.25,-1.9) {\textcolor{black}{\large (d)}};
    \end{tikzpicture}
\caption{Merging neighborhoods (highlighted in green) and overlap counts of cells in the neighborhood of a corner geometry.  
a) $M_{i+1,j-1}$ in green, and $N_{i+1,j-1} = 1$ because $(i+1,j-1)$ is not in any of the other neighborhoods; 
b) $M_{i,j-1}$   in green, and $N_{i,j-1} = 1$ because $(i,j-1)$ is not in any of the other neighborhoods; 
c) $M_{i,j}$     in green, and $N_{i,j} = 3$ because $(i,j)$ is in the neighborhoods shown in (b) and (d);
d) $M_{i+1,j}$    in green, and $N_{i+1,j} = 2$ because $(i+1,j)$ is in the neighborhood of $(i+1,j-1)$ as shown in (a).
}\label{fig:cutsWithCounts}
\end{figure}

Before the time stepping portion of the finite volume solver begins, the mesh is preprocessed by associating a merging neighborhood with each cell in the base grid, both whole and cut.  Small cut cells are merged with their neighbors until the volume of the neighborhood is at least 
$$V_{target} = \Delta x \Delta y \Delta z/2,$$
a threshold informed by results in \cite{MJB}. 
This is illustrated in two space dimensions in Figure \ref{fig:cutsWithCounts}a, where cell $(i+1,j-1)$ merges with cell $(i+1,j)$ to form  neighborhood $M_{i+1,j-1}$, highlighted in green.  
Similarly, cell $(i,j-1)$ merges with cell $(i,j)$ to form $M_{i,j-1}$'s neighborhood, shown in Figure \ref{fig:cutsWithCounts}b.
A cell with a volume larger than the threshold value does not merge with any neighboring cells, thus its merging neighborhood only contains itself, like cell $(i,j)$ in Figure \ref{fig:cutsWithCounts}c. Cell $(i+1,j)$ merges with cell $(i,j)$, shown in  Figure \ref{fig:cutsWithCounts}d.

The number of neighborhoods that contain cell $(i,j,k)$, called the overlap count $N_{i,j,k}$, is also indicated in Figure \ref{fig:cutsWithCounts} in two space dimensions.
Cell $(i,j)$ is contained in two other neighborhoods, that of cell $(i,j-1)$ and $(i+1,j)$, plus it is in its own neighborhood by definition, so its overlap count is three (Figure \ref{fig:cutsWithCounts}c).

There are several possible ways to choose which neighbors a cut cell merges with.  In this example, a cell merges in the direction closest to the boundary normal.  In the AMReX-Hydro implementation, if merging with the one cell in the direction closest to the normal does not result in a neighborhood with a volume greater than $V_{target}$, the neighbor in the direction of the next largest component of the normal is added to the neighborhood.  
In order to not create an ``L-shaped" neighborhood we then automatically add the cell in the same plane as the two neighbors that defines the neighborhood as a 2 $\times$ 2 box.  In 3D if this neighborhood is still not large enough, we then add the remaining four cells to have a 2 $\times$ 2 $\times$ 2 box. 
This is done for aesthetic reasons; an L-shaped domain would still be stable. 
This is illustrated in Figure \ref{fig:3x3}d, where normal merging for cell 9 is insufficient and a larger neighborhood must be used.
This neighborhood is useful in cases where it is desireable to preserve symmetries in the numerical solution, e.g., when the embedded boundary forms exactly a $45^{\circ}$ angle with the background Cartesian grid.

Another alternative is called central merging, where a cut cell is merged with all flow cells that lie on a $3 \times 3 \times 3$ tile centered on that cut cell.  This is neighborhood is larger and more diffusive than normal merging or the 2 $\times$ 2 $\times$ 2 box discussed above.  We will use central merging in some two-dimensional numerical examples (section \ref{sec:ex_2d}).

In AMReX-Hydro we also made a design decision to preserve spatial symmetries whenever possible. For example, if the boundary normal in a cut cell points has identical components with $n_x = n_y > n_z,$  we create a neighborhood that is larger than strictly necessary on volumetric considerations alone, to avoid directional bias in a possibly otherwise symmetric geometry and solution. 
In the above case the neighborhood would include all fluid cells in a 2 $\times$ 2 $\times$ 1 region in the cut cell's merging neighborhood.  If $n_x = n_y = n_z$ then the neighborhood would include all fluid cells in a 2 $\times$ 2 $\times$ 2 region in the cut cell's merging neighborhood.

 Once the neighborhoods have been identified, we define a neighborhood's weighted volume $\widehat V_{i,j,k}$ using the standard cell volume $V_{i,j,k}$ by
 \begin{equation}
\label{eqn:voldef}
{\widehat V}_{i,j,k} =  \sum_{(r,s,t) \in M_{i,j,k} } \,  \frac{V_{r,s,t}}{N_{r,s,t}},
\end{equation}
and a weighted centroid 
\begin{equation}
\label{eqn:centroiddef}
({\widehat x}_{i,j,k},{\widehat y}_{i,j,k},{\widehat z}_{i,j,k}) = \frac{1}{\widehat V_{i,j,k}} \sum_{(r,s,t) \in M_{i,j,k} } \,  \frac{V_{r,s,t}}{N_{r,s,t}}(x_{r,s,t},y_{r,s,t},z_{r,s,t}),
\end{equation}
where $(x_{i,j,k}, y_{i,j,k},z_{i,j,k})$ is the standard cell centroid, and $M_{i,j,k}$ is the set of cells $(r,s,t)$ in the neighborhood of cell $(i,j,k)$.  In the above, cell volumes are weighted by the inverse of their overlap counts.
A cell's overlap count, as well as the weighted volume and centroid associated to neighborhoods are required during the state redistribution postprocessing step, which we describe below.

 \subsubsection{Postprocessing}\label{sec:postprocessing}
State redistribution is implemented as a postprocessing step acting on cell averages updated by a finite volume scheme. (The finite volume scheme itself needs to be modified for the presence of cut cells to preserve accuracy, but not stability, since that is handled by SRD.) Thus, the first step of the algorithm is to compute a provisional, but possibly unstable cell update
\begin{equation}
\widehat U_{i,j,k} = U^n_{i,j,k} - \frac{\Delta t}{V_{i,j,k}}\sum_{\ell \in \text{faces}}  \mathbf{F}_{\ell}^*  \cdot \mathbf{n}_\ell A_\ell,
\label{eq:eb_update}
\end{equation}
where $A_\ell$ is a face area, $\mathbf{F}_\ell^*$ is the numerical flux, and $\mathbf n_\ell$ is the outward pointing normal.
This update is computed on all cells using a time step $\Delta t$ that is proportional to the size of cells in the background Cartesian grid. 

 The next step is to compute a weighted solution average on each merging neighborhood
 \begin{equation}
 \label{eqn:qhat}
\widehat{Q}_{i,j,k} =  \frac{1}{{\widehat V}_{i,j,k}} \, \sum_{(r,s,t) \in M_{i,j,k}} \,  
\frac{V_{r,s,t}}{N_{r,s,t}}  \widehat{U}_{r,s,t}.
 \end{equation}
For second order accuracy, a gradient is reconstructed in all cells $(i,j,k)$ for which  $M_{i,j,k}$ has at least two cells ($(i,j,k)$ itself and at least one other).
This results in a linear function for neighborhood $(i,j,k)$ of the form
 \begin{equation}\label{eq:qrecon}
\widehat{q}_{i,j,k}(x,y,z) = \widehat{Q}_{i, j,k} + \widehat{\sigma}_{x,i,j,k}(x - \widehat{x}_{i,j,k}) + \widehat{\sigma}_{y,i,j,k}(y - \widehat{y}_{i,j,k}) + \widehat{\sigma}_{z,i,j,k}(z - \widehat{z}_{i,j,k}),
\end{equation}
where $\widehat \sigma$ is the gradient obtained with a least squares approach. 
For this step the least squares equations are  centered on the weighted neighborhood centroids $(\widehat{x}_{i,j,k},\widehat{y}_{i,j,k},\widehat{z}_{i,j,k})$, not the original cell centroids
$(x_{i,j,k},y_{i,j,k},z_{i,j,k})$.

In each coordinate direction, we compute the maximum distance between the weighted centroids of the neighborhoods in the reconstruction stencil to the weighted centroid of the neighborhood on which the reconstruction is centered.
When the reconstruction stencil does not include points further than a threshold in any one coordinate direction, we grow the extent of the stencil in that coordinate direction as described in \cite{BG}.
The threshold distances are defined as $\Delta x/2$, $\Delta y/2$, $\Delta z/2$, in the $x$, $y$, and $z$ directions, respectively.
For example, if the reconstruction stencil is initially the $3\times 3 \times 3$ block, and the maximum distance in the $x$-direction $\max_{m \ne i} |\widehat{x}_i-\widehat{x}_m|$, for cells $m$ in the reconstruction stencil, is not larger than $\Delta x/2$, then the stencil is increased to $5\times 3 \times 3$.
We also point out that the $3 \times 3 \times 3$ block may not contain 27 cells as it could contain covered cells that are not included in the stencil. If the larger stencil still does not contain not enough cells, a well-conditioned gradient reconstruction is not possible, and we drop to first order.

We have also experimented with the less restrictive criterion of requiring the difference between the maximum and minimum centroid values in the stencil to be greater than half the mesh width. More experience in complicated geometries in three dimensions are necessary to ensure this is sufficiently stable.
If necessary, a Barth-Jespersen-style limiting is performed in the neighborhood gradient reconstruction step as well.

The final step is to use the neighborhood polynomial to compute the final solution update
 \begin{equation} \label{eqn:final_update_linear2}
        U^{n+1}_{i,j,k} =   \frac{1}{N_{i,j,k}}\sum_{(r,s,t)  \in W_{i,j,k}}\widehat{q}_{r,s,t}(x_{i,j,k},y_{i,j,k},z_{i,j,k}).
\end{equation}
where $W_{i,j,k}$ is the set of neighborhood indices that include cell $(i,j,k)$ in their neighborhood. 
In words, the stabilized neighborhood solution replaces the unstable finite volume update in a conservative and accurate manner.
In the above formula, the stabilized solution average on a cell in the base grid is obtained by averaging the centroid values of the overlapping merging neighborhood polynomials.  As in the two-dimensional algorithm, SRD can be used with either a method of lines integration in which the SRD algorithm is applied at each stage or with second-order Godunov type approach that directly computes fluxes at the $n+ \frac{1}{2}$ level.

This postprocessing procedure is also applied to the initial data, before any time steps have been taken. This is akin to pre-merging when using cell merging.
It is also analogous to the procedure in incompressible flow, where the initial conditions have an initial projection before time stepping begins.  We have found that our test cases have better monotonicity properties with this additional step.

\begin{figure}[ht!]
    \centering
    \includegraphics[height=.5in]{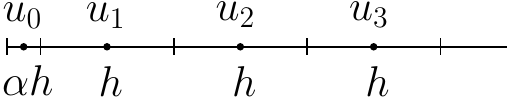}
    \caption{Notation for 1D model problem with one cut cell at the boundary.}
    \label{fig:1dBndryFig}
\end{figure}

{\bf Remark:}  This procedure can also be applied to a diffusive flux written in finite volume form. Consider the heat equation $u_t = u_{xx}. $ For this one-dimensional model problem we put one small cut cell at the left boundary of size $\alpha h$ in an otherwise uniform grid with mesh width $h$, as shown in  Figure \ref{fig:1dBndryFig}. 
One possible (though not very accurate) finite volume discretization is
\begin{equation}
u_0^{n+1} = u_0^n + \frac{\Delta t}{\alpha h} \left ( \frac{(u_1^n-u_0^n)}{(1+\alpha) h/2} - f_B) \right ) . \label{eq:diff}
\end{equation}
where $f_B$ is a specified boundary flux.
Note that it is $\alpha (1+\alpha)h^2/2 $ and not $(\alpha h)^2$ that appears in the denominator. 
Hence the same SRD procedure that eliminates  $\alpha$ from the denominator for advective updates can be used here too.  As $\alpha \rightarrow 0$, the coefficients of the stabilized finite volume scheme remain bounded. Numerical experiments for this model problem confirm this. A 3D example including a diffusive flux will be presented in sections \ref{sec:compreact} and \ref{sec:compResults}.

\section{Generalizing State Redistribution }\label{sec:alphabeta}
We present here an improvement to the original algorithm, as just described, which addresses several issues. 
The original algorithm has a sharp cutoff when the volume fraction of a cut cell reaches 0.5 and the stencil abruptly shifts. 
A transition that gracefully shuts off the amount of redistribution would be smoother and introduce less numerical dissipation. It also opens the door to consideration of weightings for other purposes (section \ref{sec:framework}).
The second issue is that the preprocessing step described earlier might generate neighborhood volumes that are overly large, resulting in more redistribution than is necessary for stability.    The modification described below addresses both these issues. The computational results in section \ref{sec:compResults} use this modification.

\subsection{A Weighted SRD Algorithm} \label{sec:framework}
The main idea is that a large cut cell that is just below the threshold volume requires less stabilization from its neighbors.  We use the notation $M^-_{i,j,k} = M_{i,j,k} - \{(i,j,k)\}$, i.e. $M^-_{i,j,k}$ is the set of cells contained in the neighborhood of cell $(i,j,k)$ with the exception of cell $(i,j,k)$ itself. As before, $W_{i,j,k}$ is the set of
$N_{i,j,k}$ indices $(r,s,t)$ containing cell $(i,j,k)$ in their neighborhood, i.e. so that cell $(i,j,k)$ is in $M_{r,s,t}.$ We again set $W_{i,j,k}^- = W_{i,j,k} - \{(i,j,k)\}$. 

For the weighted version of SRD, we define two new scalars for each cell, $\alpha_{i,j,k}$ and $\beta_{i,j,k},$ which are used as weights for the relative contributions of cell $(i,j,k)$ and its neighbors to the solution at $(i,j,k).$  
In cells with volume $V_{i,j,k} < V_{target}$ we define
\begin{equation}
{\beta}_{i,j,k} =  (V_{target} - V_{i,j,k}) \; \; /  \sum_{(r,s,t) \in M^-_{i,j,k} } \, V_{r,s,t} \, ;
\label{eqn:beta}
\end{equation}
and then 
\begin{equation}
\label{eqn:alpha}
{\alpha}_{i,j,k} =  1 \, - \,  \frac{1}{N_{i,j,k}} \, \sum_{(r,s,t) \in W^-_{i,j,k} } \, \beta_{r,s,t} \, .
\end{equation}
The idea of $\beta$ is to control the contribution a cut cell gets from its merging neighborhood. 
If a cut cell has volume $V_{i,j,k} \geq \Delta x \Delta y  \Delta z /2$, we set $\beta_{i,j,k} = 0$.
Equations \eqref{eqn:beta} and \eqref{eqn:alpha} imply  that  $0 \le \alpha_{i,j,k}, \beta_{i,j,k} \le 1$.  

As $V_{i,j,k}$
becomes larger, the weighted algorithm increases the dependence of the solution in cell $(i,j,k)$ on $\widehat{U}_{i,j,k}$ relative to the original SRD algorithm.
Notice that this version is not equivalent to the original SRD algorithm and may have some different stability properties, especially for one-dimensional test cases that are not representative of Cartesian cut cell meshes in higher dimensions.
 
All the preprocessing formulas are now modified to  pull the cut cell out of the expressions and weight it by $\alpha$, with the rest weighted by $\beta$.  Instead of Equation~(\ref{eqn:voldef}), we define a neighborhood's weighted volume $\widehat V_{i,j,k}$ using the standard cell volume $V_{i,j,k}$ by
 \begin{equation}
 \label{eqn:nbhd_alpha}
{\widehat V}_{i,j,k} =  \alpha_{i,j,k} V_{i,j,k} +  \beta_{i,j,k} \sum_{(r,s,t) \in M^-_{i,j,k}} \,  \frac{V_{r,s,t}}{N_{r,s,t}}.
\end{equation}
Instead of Equation~(\ref{eqn:centroiddef}) we use 
\begin{align}
\begin{aligned}
({\widehat x}_{i,j,k},{\widehat y}_{i,j,k},{\widehat z}_{i,j,k}) =  \frac{1}{\widehat V_{i,j,k}} \biggl(
&\alpha_{i,j,k} V_{i,j,k} (x_{i,j,k},y_{i,j,k},z_{i,j,k}) \, + \\
& \; \beta_{i,j,k}  \sum_{(r,s,t) \in M^-_{i,j,k} }  \frac{V_{r,s,t}}{N_{r,s,t}}(x_{r,s,t},y_{r,s,t},z_{r,s,t}) \biggr)
\end{aligned}
\end{align}
to define the weighted centroid.
 
 In the postprocessing step, instead of Equation~(\ref{eqn:qhat}), we define the
 weighted solution average on each merging neighborhood as
 \begin{equation}
 \label{eqn:qhat_alpha}
\widehat{Q}_{i,j,k} =  \frac{1}{{\widehat V}_{i,j,k}} \, \left( \alpha_{i,j,k} V_{i,j,k} \widehat{U}_{i,j,k} \, + \beta_{i,j,k} \,
\sum_{(r,s,t) \in M^-_{i,j,k}}   
\frac{V_{r,s,t}}{N_{r,s,t}}  \widehat{U}_{r,s,t} \right).
 \end{equation}

Slopes are computed as in the original algorithm, and we use Equation~(\ref{eq:qrecon}) to construct $\widehat{q}$ as before.  The final step is to use the neighborhood polynomial to compute the final solution update. Instead of Equation~(\ref{eqn:final_update_linear2}) we use
 \begin{equation} \label{eqn:final_update_alpha}
        U^{n+1}_{i,j,k} =   \alpha_{i,j,k} \widehat{q}_{i,j,k} + 
        \frac{1}{N_{i,j,k}}\sum_{(r,s,t)  \in W^-_{i,j,k}} \beta_{r,s,t} \; 
        \widehat{q}_{r,s,t}(x_{i,j,k},y_{i,j,k},z_{i,j,k}).
\end{equation}
We will show that this new variant retains the linearity and conservation properties of the original algorithm. 
\begin{figure}
\centering
\begin{tikzpicture}
\node (1) at (0,0) {  \begin{tikzpicture}

\node(aa) at (-0.85, 4.5) {cell index};

\clip (-2.1,0.65) rectangle (0.9,4.1);
\draw (-6,-6) grid (6,6);
\path [fill=white] (-4,1.85) rectangle (5,-5);
\path [fill=white] (0.4,0) rectangle (6,6);
\draw[line width=0.4mm] (-3,1.85) -- (0, 1.85) -- (0.4,2) -- (0.4, 6); 
\node (0) at (-0.75, 1.) {(a)};

\node (a) at (-1.5,3.5) {\scriptsize $1$};
\node (b) at (-0.5,3.5) {\scriptsize $2$};
\node (c) at ( 0.7 ,3.5) {\scriptsize $3$};
\node (d) at (-1.5,2.5) {\scriptsize $4$};
\node (e) at (-0.5,2.5) {\scriptsize $5$};
\node (f) at ( 0.7 ,2.5) {\scriptsize $6$};
\node (d) at (-1.5,1.5) {\scriptsize $7$};
\node (e) at (-0.5,1.5) {\scriptsize $8$};
\node (f) at ( 0.55 ,1.5) {\scriptsize $9$};

\draw (0.45 , 1.6) -- (0.1,  1.95);
\draw (-0.5 ,1.6) -- (-0.5, 1.95);
\draw (-1.5 ,1.6) -- (-1.5, 1.95);

\draw (0.55, 3.5) -- (0.2, 3.5);
\draw (0.55, 2.5) -- (0.2, 2.5);

  \end{tikzpicture}};
\node (2) at (0.1+3.75,0) {  \begin{tikzpicture}

\node(aa) at (-0.85, 4.5) {};
\clip (-2.1,0.76) rectangle (0.9,4.1);
\draw (-6,-6) grid (6,6);
\path [fill=white] (-4,1.85) rectangle (5,-5);
\path [fill=white] (0.4,0) rectangle (6,6);
\draw[line width=0.4mm] (-3,1.85) -- (0, 1.85) -- (0.4,2) -- (0.4, 6); 
\node (0) at (-0.75, 1.) {(b)};

\draw [ultra thick, draw=green!70!black, fill=green!70!black, fill opacity=0.2] (0, 1.85) rectangle  (-1,3);

\node (e) at (-0.5,1.5) {\scriptsize $8$};

\draw (-0.5 ,1.6) -- (-0.5, 1.95);

  \end{tikzpicture}};
\node (3) at (0.1+7.5,0) {  \begin{tikzpicture}

\node(aa) at (-0.85, 4.5) {};
\clip (-2.1,0.76) rectangle (0.9,4.1);
\draw (-6,-6) grid (6,6);
\path [fill=white] (-4,1.85) rectangle (5,-5);
\path [fill=white] (0.4,0) rectangle (6,6);
\draw[line width=0.4mm] (-3,1.85) -- (0, 1.85) -- (0.4,2) -- (0.4, 6); 
\node (0) at (-0.75, 1.) {(c)};

\draw [ultra thick, draw=green!70!black, fill=green!70!black, fill opacity=0.2] (0.4,2) rectangle  (-1,3);

\node (f) at ( 0.7 ,2.5) {\scriptsize $6$};

\draw (0.55, 2.5) -- (0.2, 2.5);

  \end{tikzpicture}};
\node (4) at (0.1+11.25,0) {  \begin{tikzpicture}

\node(aa) at (-0.85, 4.5) {};
\clip (-2.1,0.76) rectangle (0.9,4.1);
\draw (-6,-6) grid (6,6);
\path [fill=white] (-4,1.85) rectangle (5,-5);
\path [fill=white] (0.4,0) rectangle (6,6);
\draw[line width=0.4mm] (-3,1.85) -- (0, 1.85) -- (0.4,2) -- (0.4, 6); 
\node (0) at (-0.75, 1.) {(d)};

\draw [ultra thick, draw=green!70!black, fill=green!70!black, fill opacity=0.2] (0.4,2) -- (0, 1.85) -- (-1, 1.85)-- (-1, 3) -- (0.4, 3) -- (0.4, 2) -- (0, 1.85);

\draw[ultra thick, draw=red!70!black, opacity=1, dashed] (0, 1.85) -- (0.4,2) -- (0.4, 3) -- (0, 3) -- (0, 1.85);

\node (f) at ( 0.55 ,1.5) {\scriptsize $9$};

\draw (0.45 , 1.6) -- (0.1,  1.95);

  \end{tikzpicture}};
\node (5) at (7.25, 2.05) {\resizebox{0.32\textwidth}{!}{\Large merging neighborhoods}};
\end{tikzpicture}

\caption{Figure (a) illustrates a configuration where cell 3 is merged with cell 2, cell 6 is merged with cell 5, cell 7 is merged with cell 4, cell 8 is merged with cell 5 and cell 9 is merged with cells 5, 6 and 8. Cell 5 is contained in 3 other neighborhoods as well as its own, resulting in an overlap count of 4.  Cell 1 is in only its own neighborhood.
Figures (b), (c), and (d) indicate the merging neighborhoods associated with cells 6, 8, and 9, respectively.
For cells 6 and 8, normal merging is sufficient.  However, for cell 9, normal merging results in the neighborhood indicated by the red dashed line in (d).
When normal merging is insufficient, we use the larger green neighborhood instead.
}\label{fig:3x3}
\end{figure}
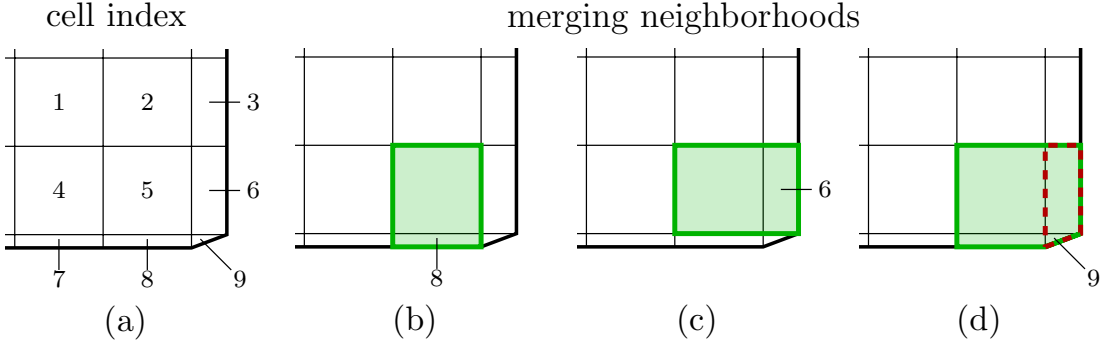

To show conservation, define the matrix $A$ such that
\[
\widehat{V} = A V,
\]
where $V$ is the vector of cell volumes and $\widehat{V}$ is the vector of neighborhood volumes. To make it concrete, for the case considered in Figure \ref{fig:3x3}, the matrices $A$ for the original and weighted algorithms are given by 

\[
\renewcommand*{\arraystretch}{1.2}
A^\mathrm{orig} = 
\begin{bmatrix}
1 & 0 & 0 & 0 & 0 & 0 & 0 & 0 & 0 \\
0 & \frac{1}{2} & 0 & 0 & 0 & 0 & 0 & 0 & 0 \\
0 & \frac{1}{2} & 1 & 0 & 0 & 0 & 0 & 0 & 0 \\
0 & 0 & 0 & \frac{1}{2} & 0 & 0 & 0 & 0 & 0 \\
0 & 0 & 0 & 0 &  \frac{1}{4} & 0 & 0 & 0 & 0 \\
0 & 0 & 0 & 0 &  \frac{1}{4} & \frac{1}{2} & 0 & 0 & 0 \\
0 & 0 & 0 & \frac{1}{2} & 0 & 0 & 1 & 0 & 0 \\
0 & 0 & 0 & 0 &  \frac{1}{4} & 0 & 0 & \frac{1}{2} & 0 \\
0 & 0 & 0 & 0 &  \frac{1}{4} & \frac{1}{2} & 0 & \frac{1}{2} & 1 
\end{bmatrix}
\quad \text{and} \quad
A^\mathrm{wght} = 
\begin{bmatrix}
1 & 0 & 0 & 0 & 0 & 0 & 0 & 0 & 0 \\
0 & {\alpha_2} & 0 & 0 & 0 & 0 & 0 & 0 & 0 \\
0 & \frac{\beta_3}{2} & 1 & 0 & 0 & 0 & 0 & 0 & 0 \\
0 & 0 & 0 & {\alpha_4} & 0 & 0 & 0 & 0 & 0 \\
0 & 0 & 0 & 0 &  {\alpha_5} & 0 & 0 & 0 & 0 \\
0 & 0 & 0 & 0 &  \frac{\beta_6}{4} & {\alpha_6} & 0 & 0 & 0 \\
0 & 0 & 0 & \frac{\beta_7}{2} & 0 & 0 & 1 & 0 & 0 \\
0 & 0 & 0 & 0 &  \frac{\beta_8}{4} & 0 & 0 & {\alpha_8} & 0 \\
0 & 0 & 0 & 0 &  \frac{\beta_9}{4} & \frac{\beta_9}{2} & 0 & \frac{\beta_9}{2} & 1 
\end{bmatrix}.
\]

In $A^{\mathrm{orig}}$ the entries $a_{ij}$ are the inverse of cell $j$'s overlap count, if cell $i$ is in cell $j$'s neighborhood. In $A^\mathrm{{wght}}$, the generalized weights correspond to the previously defined $\alpha$ and $\beta$.
The key property of both matrices  $A$ is that the columns sum to 1, 
\[
e^T A = e^T,
\]
where $e$ is a vector of all 1's. 
This is the property that any state redistribution method must satisfy to be conservative.
We then have that
\[
\widehat{Q} = \mathrm{Diag}(\widehat{V})^{-1} A \;\mathrm{Diag}(V) \widehat{U}.
\]
Furthermore, if slopes are zeroed in Eq. (\ref{eq:qrecon}) then
\[
U^{n+1} = A^T \widehat{Q}.
\]
Conservation then follows from
\begin{align*}
V^T U^{n+1} &= V^T A^T \; \mathrm{Diag}(\widehat{V})^{-1} A \; \mathrm{Diag}(V) \widehat {U}
= 
 \widehat{V}^T \; \mathrm{Diag}(\widehat{V})^{-1} A \; \mathrm{Diag}(V) \widehat {U}, \\
 &= e^T A \; \mathrm{Diag}(V) \widehat{U} = V^T \widehat{U} .
\end{align*}

Including the slopes in Eq. (\ref{eq:qrecon}) redistributes mass to preserve linearity but does not alter the conservation properties of the method. Continuing, 
in matrix form using \eqref{eq:qrecon} the update with slopes can be written as
\begin{align*}
U^{n+1} = A^T Q + \left[\mathrm{Diag}(x)A^T-A^T \mathrm{Diag}(\widehat{x}) \right ]\sigma_x
 &+\left[ \mathrm{Diag}(y)A^T-A^T \mathrm{Diag}(\widehat{y}) \right ]\sigma_y \\
 &+ \left[ \mathrm{Diag}(z)A^T-A^T \mathrm{Diag}(\widehat{z}) \right ]\sigma_z,
\end{align*}
where $x,y,z$ are the centroids of the original cells, $\widehat{x},\widehat{y},\widehat{z}$
are the weighted centroids of the neighborhoods and $\sigma_x,\sigma_y, \sigma_x$ are the slopes in the neighborhood reconstruction.
We then note that
\begin{align*}
V^T \left[\mathrm{Diag}(x)A^T-A^T \mathrm{Diag}(\widehat{x}) \right ]\sigma_x  &=
\left [ x^T \mathrm{Diag}(V)A^T- V^T A^T \mathrm{Diag}(\widehat{x}) \right ]\sigma_x \\
&= \left [ \widehat{x}^T \mathrm{Diag}(\widehat{V}) - \widehat{V}^T \mathrm{Diag}(\widehat{x}) \right ] \sigma_x = 0.
\end{align*}

Although not pursued here, the structure suggests that any collection of weights with
$e^T A = e^T$ can be used to define the SRD algorithm.  This opens the possibility of more sophisticated approaches, including ones that might depend on the local solution or maintain monotonicity.

\subsection{Two-dimensional numerical examples}\label{sec:ex_2d}
We present two two-dimensional  examples showing the performance of  weighted SRD.  First, we solve the linear advection equation and show the superior performance of both the original and weighted SRD algorithms relative to flux redistribution, which was the stabilization technique in AMReX.
We also show that the weighted algorithm smoothly activates as the volume fraction decreases below a threshold value.
For a discontinuous solutions, this will prevent $\mathcal{O}(1)$ changes in the numerical solution under $\mathcal{O}(\epsilon)$ changes to the mesh.
Second, we present a convergence study using an analytic solution of the Euler equations, showing that the error can be reduced by a factor of between 2 and 3 for central merging when using the weighted algorithm.  We even observe a 10 to 15\% reduction in the error when using the less diffusive normal merging neighborhood.
Finally, we show that both original and weighted algorithms converge at the same rate as the base finite volume scheme.

\subsubsection*{Linear advection} \label{sec:linadv2d}
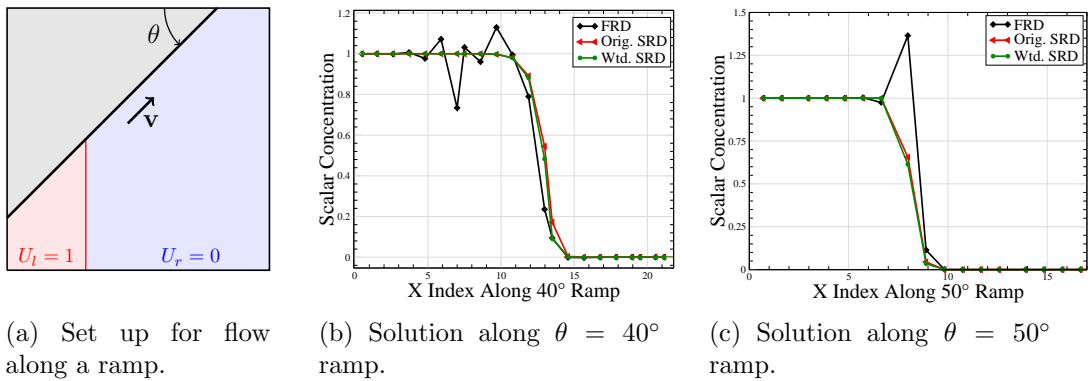
\begin{figure}[ht]
    \centering
    \begin{subfigure}[t]{.24\textwidth}
     \centering
       \begin{tikzpicture}

\draw[line width=0.5mm, fill, red!10!white] (0,0) -- ++ (1.5, 0) -- ++(0,2.5) --++(-1.5, -1.5);
\draw[line width=0.5mm, fill, blue!10!white] (1.5, 0) -- ++(3.5, 0) --++(0,5) --++(-1,0) --++(-2.5, -2.5);
\draw[line width=0.25mm, solid, red] (1.5,0.0) --++(0,2.5);
\node (AA) at (0.75, 0.25) {\textcolor{red}{$U_{l}=1$}};
\node (BB) at (3.5, 0.25) {\textcolor{blue}{$U_{r}=0$}};

\draw[line width=0.5mm, fill, black!10!white] (0,1) -- ++ (4,4) -- ++(0,0) -- ++(-4,0) --++(0,-4);
\draw[line width=0.3mm] (0, 0) -- (5, 0) -- (5,5) -- (0, 5) -- (0,0);

\draw[line width=0.5mm] (0,1) -- ++(4,4);

\begin{scope}[shift={(1.3,0.3)}]
\draw[line width=0.5mm, ->] (1,2.5) -- ++(0.5, 0.5) ;
\node at (1.1+0.35,2.9-0.35) {\Large $\mathbf{v}$};
\end{scope}

\coordinate (A) at (3,4);
\coordinate (B) at (4,5);
\coordinate (C) at (3,5);
\draw pic[draw, "\Large $\theta$", ->, angle eccentricity=1.3, angle radius = 1cm, line width=0.25mm] {angle = C--B--A};

\node (0) at (0.1,-0.5) {};
  \end{tikzpicture}
     \caption{Set up for flow along a ramp.}\label{fig:ramp4050a}
    \end{subfigure}
    \hspace*{.15in}
    \begin{subfigure}[t]{.31\textwidth}
     \centering
     \includegraphics[height=1.58in,trim= 12 0 0 0,clip]{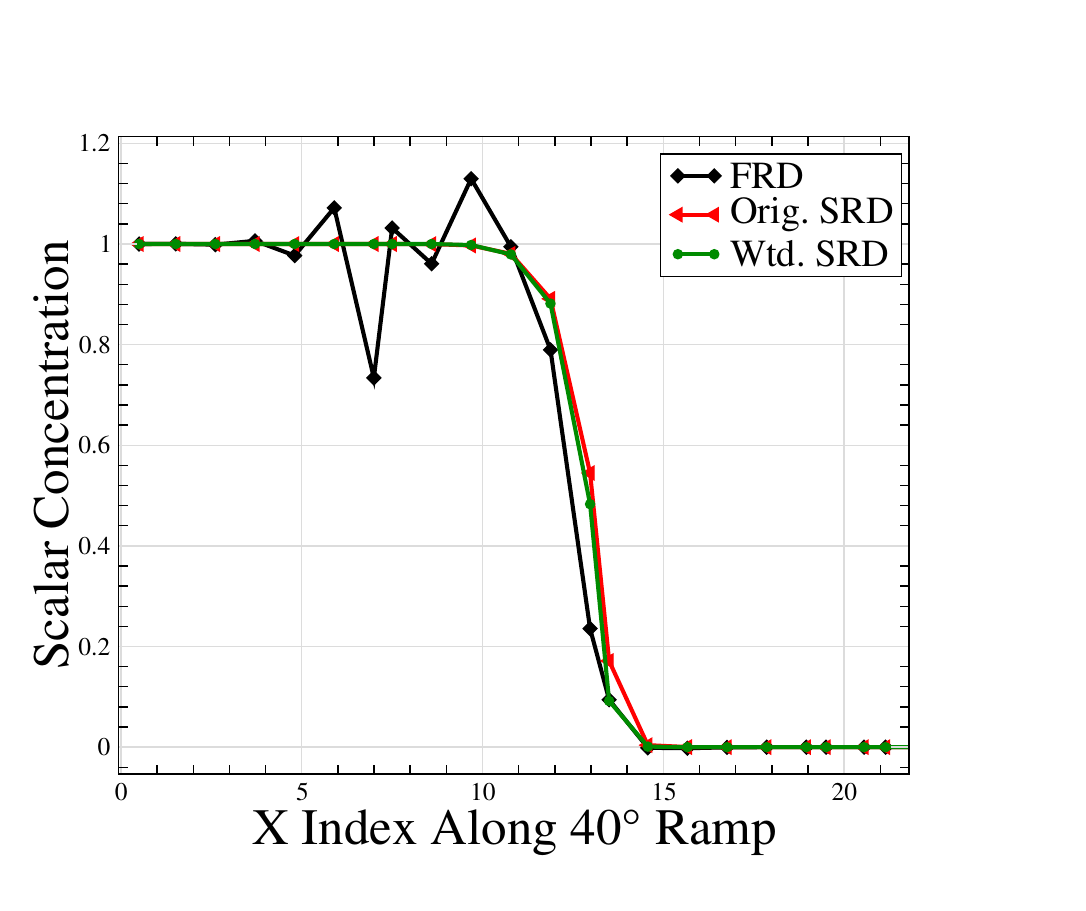}
     \caption{Solution along $\theta=40^{\circ}$ ramp.}\label{fig:ramp4050b}
    \end{subfigure}
    \hspace*{.15in}
    \begin{subfigure}[t]{.31\textwidth}
     \centering
     \includegraphics[height=1.58in,trim=12 25 50 50,clip]{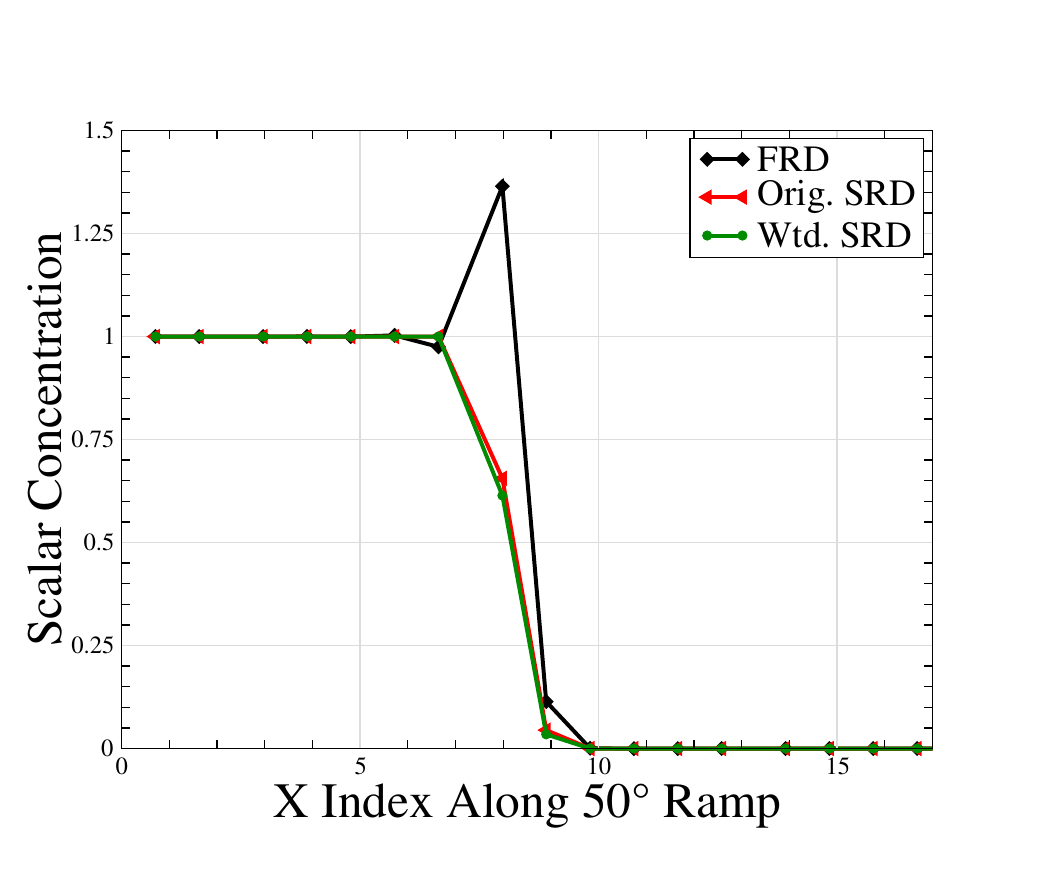}
     \caption{Solution along $\theta=50^{\circ}$ ramp.}\label{fig:ramp4050c}
     \end{subfigure}
     \caption{Linear advection of a passive scalar along a slanted wall.   FRD shows undershoots and overshoots, whereas both SRD algorithms do not. The profiles are taken after 10 time steps.
     }\label{fig:ramp4050}
\end{figure}

In this example, we consider linear advection of a passive scalar in the presence of a ramp (Figure \ref{fig:ramp4050a}), with cut cells where they intersect the ramp.  We consider several different angles of tilt ($\theta = 40$, 45 and 50 degrees) and different combinations of volume fractions.  The velocity vector $\mathbf v$ in each numerical test is parallel to the ramp.
The initial condition for the advected scalar is a discontinuous Heaviside function; the scalar has concentration 1 to the left of the discontinuity and 0 to the right. 
The flux through the embedded boundary is zero since the advection velocity is parallel to the wall.  
The values of 1 and 0 are applied to ghost cells on the left (inflow) and right (outflow) boundaries respectively.
We use constant extrapolation along the bottom boundary.

In Figures \ref{fig:ramp4050b}, \ref{fig:ramp4050c},  the solution in the cut cells is reconstructed linearly to the centroid of the embedded boundary face in each cut cell. 
Figure~\ref{fig:ramp4050b} shows the solution on a 40 degree ramp, so that the merging neighborhood is in the vertical direction. Figure \ref{fig:ramp4050c} shows the solution for a 50 degree ramp, so that the merging neighborhood is in the horizontal direction. In both examples the results after 10 time steps are shown. None of the algorithms guarantee monotonicity.  
We see however that in both the 40 degree and 50 degree cases that FRD  generates  overshoots and undershoots while both the original and weighted SRD algorithm do not.  The new weighted algorithm is seen to be less dissipative than the original.

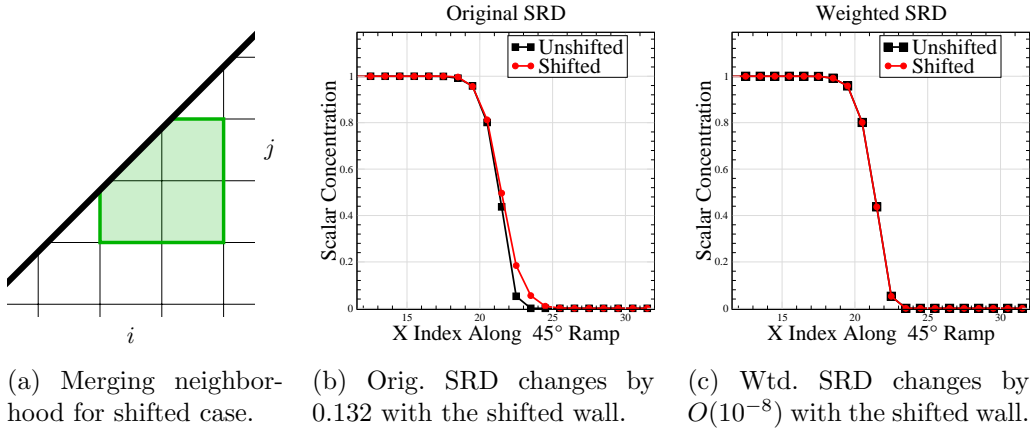
\begin{figure}[ht!]
    \centering
    \begin{subfigure}[t]{.25\textwidth}
     \centering
       \begin{tikzpicture}

\node (b) at (2.5,-0.5) {$i$};
\node (a) at (4.75,2.5) { $j$};

\clip (0.5,-0.2) rectangle (4.5,4.5);
\draw (-6,-6) grid (5,5);
\draw [ultra thick, draw=green!70!black, fill=green!70!black, fill opacity=0.2] (2,3) rectangle  (4,1);

\begin{scope}[rotate=45]
\path [fill=white] (0,-0.1) rectangle (10,5);
\end{scope}

\draw[line width=1mm] (0,-0.15) -- ++(5,5);
\draw[line width=1mm] (0,-0.15) -- ++(-1,-1);

  \end{tikzpicture}
     \caption{Merging neighborhood for shifted case.}\label{fig:eps_up}
    \end{subfigure}
    \hspace*{.05in}
    \begin{subfigure}[t]{.31\textwidth}
     \centering
     \includegraphics[width=\linewidth,trim=0 20 48 0,clip]{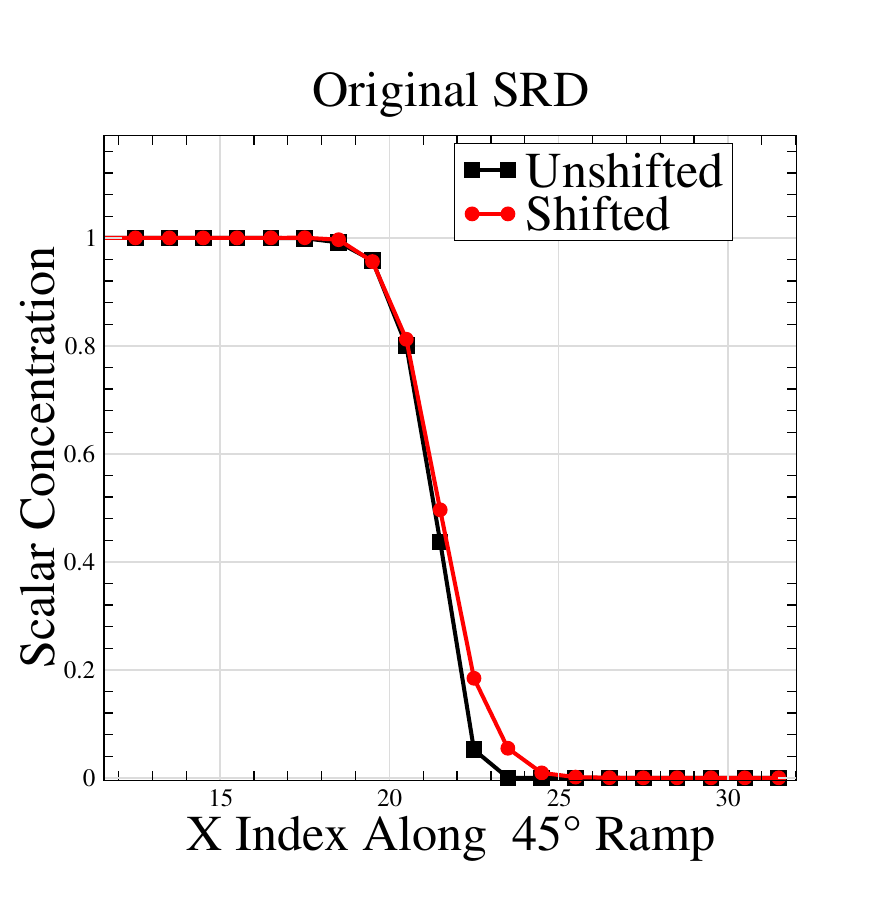}
     \caption{Orig.~SRD changes by 0.132  with the shifted wall.}
     \label{fig:45_orig}
    \end{subfigure}
    \hspace*{.07in}
    \begin{subfigure}[t]{.31\textwidth}
     \centering
     \includegraphics[width=\linewidth,trim=0 20 48 0,clip]{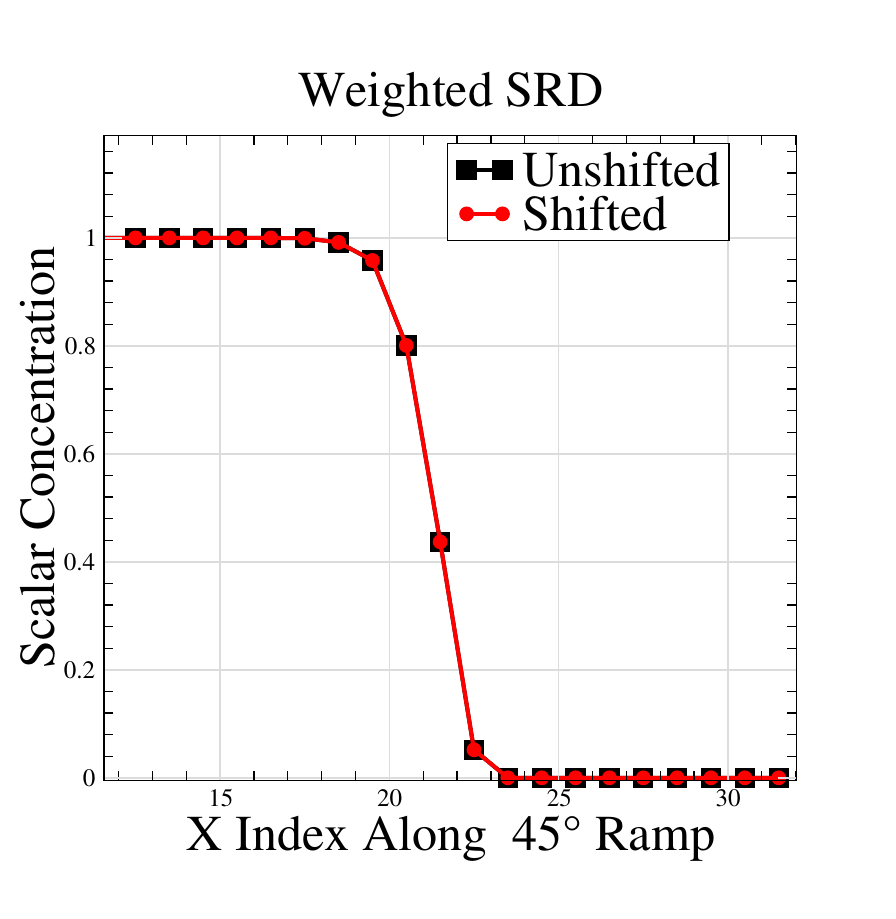}
     \caption{Wtd.~SRD changes by $O(10^{-8})$  with the shifted wall.}
      \label{fig:45_wtd}
     \end{subfigure}
     \caption{Change in solution with $2 \times 2$ merging when a $45^\circ$ ramp is shifted vertically down  by $O(10^{-7})$.
     In (b), we observe that the original SRD solution is diffused in the shifted configuration.
     In (c), no additional diffusion is introduced by the stabilization.
     }
     \label{fig:ramp45}
\end{figure}

In the next linear advection test  we compare the original and weighted SRD algorithm for same discontinuous flow problem as above, but for a ramp at exactly 45 degrees from the horizontal. If the ramp exactly bisects the Cartesian cells (assuming equal mesh widths), the cut cell volume fraction is exactly 0.5, so neither the original or weighted SRD algorithms are invoked.  
If the wall is shifted down by $O(10^{-7}) \Delta y$, the volume fractions of cut cells are either slightly below 0.5 or almost 1. 
The $2 \times 2$ merging neighborhood for this case is shown in Figure \ref{fig:eps_up}.
Recall that the
original algorithm has a sharp cutoff when the volume fraction of a cut cell decreased below 0.5 and abruptly activates, while the weighted version activates gradually.   In Figure \ref{fig:45_orig} the maximum difference between the solutions is 0.132, whereas in Figure \ref{fig:45_wtd}  
it is $O(10^{-8})$, and the 
unshifted and shifted solutions are visually indistinguishable. 
Even with smaller merging neighborhoods, e.g. if the original algorithm merges only in the vertical direction or only in the horizontal direction, the maximum difference in the solution is respectively $0.0051$, and $0.1165$.  This is still orders of magnitude larger than the $\mathcal{O}(10^{-8})$ change in the solution from the weighted algorithm.

\subsubsection*{Supersonic vortex}
\begin{figure}[ht!]
    \centering
    \begin{subfigure}[t]{.420\textwidth}
    \includegraphics[height=2.2in]{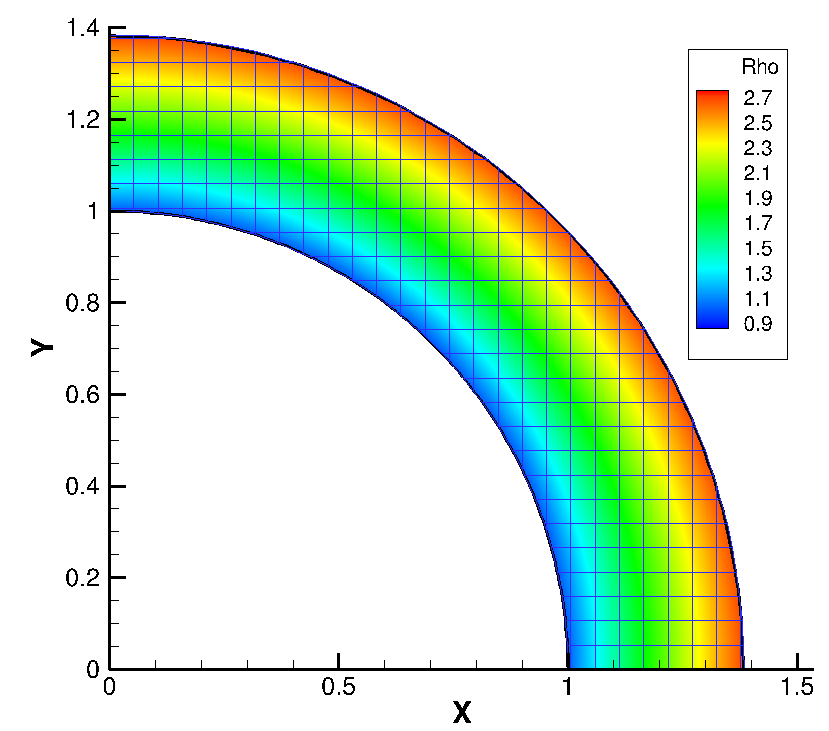} 
    \caption{Computational mesh using a $27\times27$ background grid and plot of density $\rho$ in exact solution.}
    \end{subfigure}
    \hfill
    \begin{subfigure}[t]{.5\textwidth}
    \includegraphics[height=2.2in]{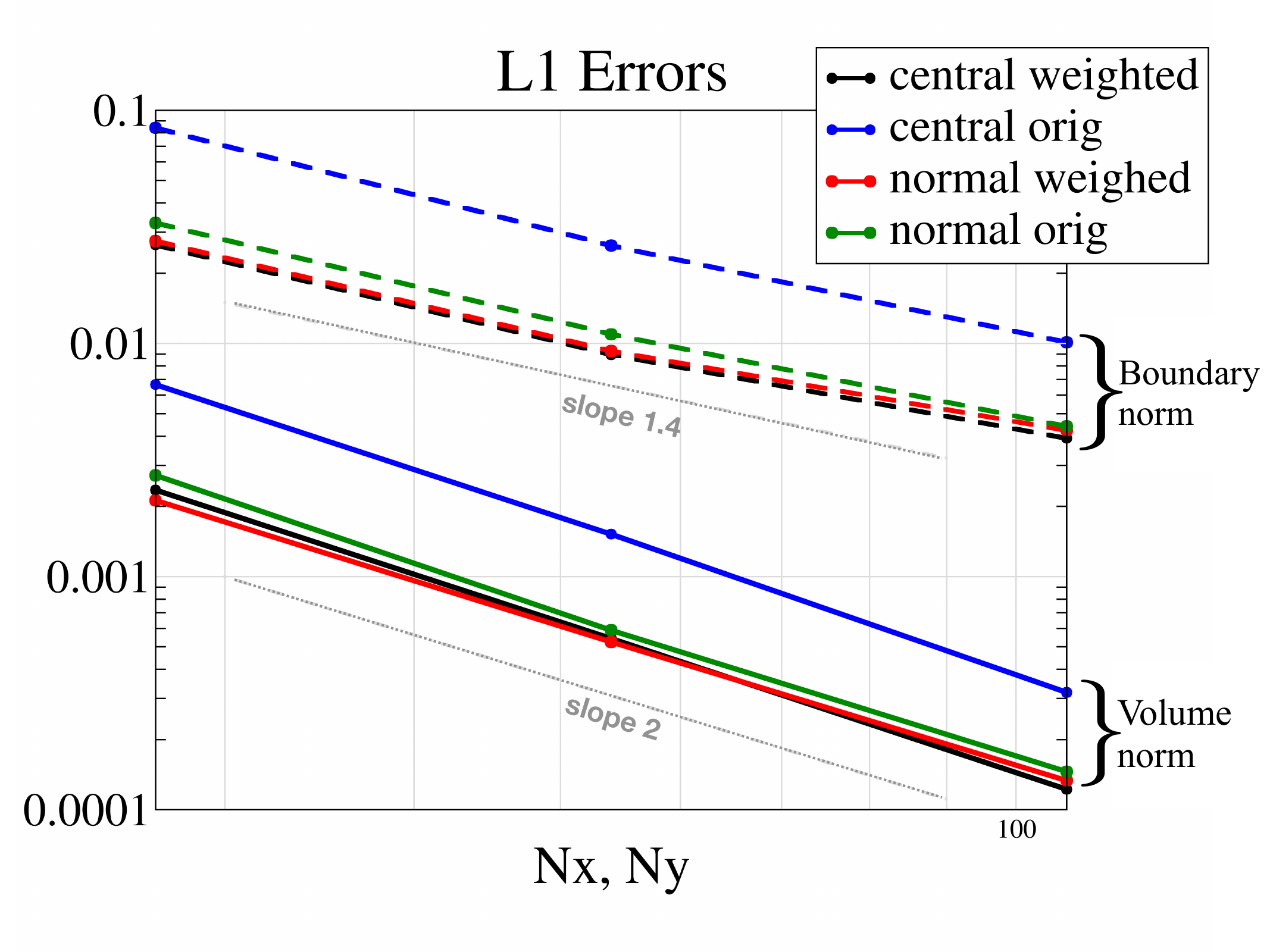}
    \caption{$L_1$ errors in density for both weighted and original algorithms on the volume and along boundary using both central and normal merging neighborhoods.}\label{fig:errors}
    \end{subfigure}
    \caption{Supersonic vortex test problem set up and convergence study. 
    For central merging, the weighted algorithm reduces the error by a factor between 2 and 3 relative to the original algorithm.
    Even for the less diffusive normal merging neighborhoods, the error is improved by 10-15\%.
    The weighted algorithm is much less sensitive to the size of the merging neighborhood than the original algorithm.
    }
    \label{fig:abOrig2andOnly_2}
\end{figure}

In this example, we solve the Euler equations for supersonic flow through a quarter annulus (Figure \ref{fig:abOrig2andOnly_2} left), with inner radius, $r_i = 1.0$, and outer radius $r_o = 1.384$.
This problem has the exact solution \cite{aftosmis:acc}
\begin{equation}
\rho = \rho_i \left \{ 1 + \frac{\gamma-1}{2} \, M_i^2 \left [ 1 - \left(\frac{r_i}{r} \right)^2
\right ] \right \}^{\frac{1}{\gamma-1}}
\end{equation}
and $ u = a_i \, M_i \, (\frac{r_i}{r})\,  \sin (\theta)$, 
$ v = -a_i\,  M_i\,  \left(\frac{r_i}{r} \right)\,  \cos(\theta)$, and
$ p = \rho^\gamma / \gamma$.
We set $\gamma = 1.4$, $\rho_i=1$, and the Mach number on the inner circle
$M_i = 2.25$.
At the inflow and outflow boundaries, the exact solution is prescribed in ghost cells.
The curved inner and outer walls are reflecting, where the the flux along the boundary is set to be
$
\mathbf{F}^* =  \begin{bmatrix}
0 ,
p n_x,
p n_y ,
0
\end{bmatrix}^T.
$
Here $p$ is the pressure of the numerical solution at the boundary, and $\mathbf{n} = (n_x, n_y)$ is the outward pointing normal.
We solve this problem on a sequence of embedded boundary grids on the domain $[0, 1.43]\times[0, 1.4301]$.

Figure \ref{fig:abOrig2andOnly_2} right shows a convergence study for both the original and weighted algorithms, where we provide the $L_1$ error of density, $\rho$, measured in the volume
$\sum_{i,j} |e_{i,j}| V_{i,j}$ and on the boundary $\sum_{i,j \in \text{bdry.}} |e_{i,j}| A_{i,j}$, where $V_{i,j}$ is the cell area and $A_{i,j}$ is the boundary edge segment length.
The error, $e_{i,j}$, in the volume and boundary formula is computed by taking the difference between the solution average on $i,j$ and the exact solution at $i,j$'s centroid.

We observe that the original algorithm is highly sensitive to the neighborhood sizes.  The larger neighborhoods of central merging result in much larger errors.
The weighted algorithm avoids this problem.
The error is reduced by a factor between 2 and 3 using the weighted algorithm with central merging neighborhoods.
The weighted algorithm reduces the dissipation to a level comparable with normal merging.
We even observe a 15\% improvement in the boundary error on coarse grids and a 10\% improvement in the volume when using weighted SRD with normal merging.
On all the grids, only a small number of cut cells had volume fractions between 0.3 and 0.5, where the weighted algorithm is most effective.
In these runs, we used second order accurate gradients at the cut cells. 

In all cases shown in Figure \ref{fig:errors}, the $L_1$ norm of the solution converges with second order; the rate at the boundary drops to approximately 1.4. 
Loss of accuracy in the numerical solution at the embedded boundary has been reported before in \cite{BG, engwer2020stabilized, nemec_tm14} and is due to the irregularity of the mesh at the cut cells.

\section{AMReX Implementation Details}\label{sec:par}
The generalized SRD algorithm is implemented for two and three spatial dimensions in C++ in the
{\it AMReX}\footnote{https://github.com/AMReX-Codes/AMReX, 2021 } software framework~\cite{amrex-ijhpc}.  It is publicly available in the open source 
{\it AMReX-Hydro}\footnote{https://github.com/AMReX-Codes/AMReX-Hydro, 2021} 
collection  of routines.  AMReX-Hydro includes modules for flux redistribution as well as state redistribution; the redistribution modules are used by compressible and low Mach number flow codes and the examples given in this paper all use this implementation.  

The fundamental data structure in AMReX  contains multi-dimensional arrays on logically rectangular patches within the computational domain.   AMReX distributes patches of data to different MPI ranks; 
communication between patches  at the same refinement level occurs most frequently by filling ghost cells (also known as halo cells);  the SRD algorithm exploits the native AMReX routines for doing so.

While the SRD algorithm can be used in a simulation code that computes the 
solution on an AMR hierarchy,  the current implementation assumes that the 
coarse-fine boundary does not cross any cut cells.  Thus the SRD algorithm is 
effectively used one level at a time, and is exactly the same 
as it would be in a code with uniform mesh resolution.

To implement the SRD algorithm on a set of rectangular patches that are all at the same level, 
we must know how many ghost cells need to be filled for each type of data. The version of the SRD algorithm implemented in AMReX starts with a  $3 \times 3 \times 3$ set of neighborhoods, no neighbor is more than one cell away in any coordinate direction. If necessary, the slope computation can increase in a direction up to a $5 \times 5 \times 5$ stencil if doing so does not require additional ghost cells.
We break the SRD algorithm into two sections for this purpose.  In the preprocessing step of SRD, we compute
the necessary geometric information such as neighborhood volumes, centroid locations, number of
neighbors, etc. This information does not change between time steps unless a regridding operation 
is performed in an AMR simulation.

\begin{figure}[ht!]
\vspace*{-.13in}
\centering
\begin{tikzpicture}
\node (11) at (-2.7,-0.9) {  \begin{tikzpicture}

\clip (-2.9,0.9) rectangle (4.1,5.1);
\draw (-6,-6) grid (6,6);
\draw [ultra thick, draw=green!70!black, fill=green!70!black, fill opacity=0.2] (-2, 1.) rectangle  (1,3);

\begin{scope}[rotate=0]
\path [fill=white] (-4,1.7) rectangle (5,-5);
\end{scope}
\path [fill=white] (-5.5,1.7) rectangle (6,-6);
\begin{scope}[rotate=-16.3]
\path [fill=white] (-6,2.62) rectangle (6,6);
\end{scope}

\definecolor{beaublue}{rgb}{0.74, 0.83, 0.9};
\draw [ultra thick, draw=beaublue, fill=beaublue, fill opacity=0.2] (-1,1.7) -- (-1,1) -- (0,1) -- (0,1.7);
\node (a) at (-0.5,1.35) {\footnotesize $i+1, j$};

\draw[line width=0.5mm] (-3,1.7) -- (0, 1.7) -- (3.45,1.7)  -- (-3,3.6);

  \end{tikzpicture}};
\node (11) at (5.2,-0.9) {  \begin{tikzpicture}

\clip (-2.9,0.9) rectangle (4.1,5.1);
\draw (-6,-6) grid (6,6);
\draw [ultra thick, draw=green!70!black, fill=green!70!black, fill opacity=0.2] (-1, 1.) rectangle  (2,3);

\begin{scope}[rotate=0]
\path [fill=white] (-4,1.7) rectangle (5,-5);
\end{scope}
\path [fill=white] (-5.5,1.7) rectangle (6,-6);
\begin{scope}[rotate=-16.3]
\path [fill=white] (-6,2.62) rectangle (6,6);
\end{scope}

\definecolor{beaublue}{rgb}{0.74, 0.83, 0.9};
\draw [ultra thick, draw=beaublue, fill=beaublue, fill opacity=0.2] (0,1.7) -- (0,1) -- (1,1) -- (1,1.7);
\node (a) at (0.5,1.35) {\footnotesize $i+2, j$};

\draw[line width=0.5mm] (-3,1.7) -- (0, 1.7) -- (3.45,1.7)  -- (-3,3.6);

  \end{tikzpicture}};
\node (11) at (-2.7,-4.9) {  \begin{tikzpicture}

\clip (-2.9,0.9) rectangle (4.1,5.1);
\draw (-6,-6) grid (6,6);
\draw [ultra thick, draw=green!70!black, fill=green!70!black, fill opacity=0.2] (0, 1.) rectangle  (3,3);

\begin{scope}[rotate=0]
\path [fill=white] (-4,1.7) rectangle (5,-5);
\end{scope}
\path [fill=white] (-5.5,1.7) rectangle (6,-6);
\begin{scope}[rotate=-16.3]
\path [fill=white] (-6,2.62) rectangle (6,6);
\end{scope}

\definecolor{beaublue}{rgb}{0.74, 0.83, 0.9};
\draw [ultra thick, draw=beaublue, fill=beaublue, fill opacity=0.2] (1,1.7) -- (1,1) -- (2,1) -- (2,1.7);
\node (a) at (1.5,1.35) {\footnotesize $i+3, j$};

\draw[line width=0.5mm] (-3,1.7) -- (0, 1.7) -- (3.45,1.7)  -- (-3,3.6);

  \end{tikzpicture}};
\node (11) at (5.2,-4.9) {  \begin{tikzpicture}

\clip (-2.9,0.9) rectangle (4.1,5.1);
\draw (-6,-6) grid (6,6);

\draw [ultra thick, draw=green!70!black, fill=green!70!black, fill opacity=0.2] (1, 1.) rectangle  (4,3);

\begin{scope}[rotate=0]
\path [fill=white] (-4,1.7) rectangle (5,-5);
\end{scope}
\path [fill=white] (-5.5,1.7) rectangle (6,-6);
\begin{scope}[rotate=-16.3]
\path [fill=white] (-6,2.62) rectangle (6,6);
\end{scope}

\definecolor{beaublue}{rgb}{0.74, 0.83, 0.9};
\draw [ultra thick, draw=beaublue, fill=beaublue, fill opacity=0.2] (2,1.7) -- (2,1) -- (3,1) -- (3,1.7);
\node (a) at (2.5,1.35) {\footnotesize $i+4, j$};

\draw[line width=0.5mm] (-3,1.7) -- (0, 1.7) -- (3.45,1.7)  -- (-3,3.6);

  \end{tikzpicture}};

\node (11) at (-3,-3.35) {\large (a) \begin{tabular}{l}
$\widehat{Q}_{i+1}$ and gradient needed\\
to update cell $i$.
\end{tabular} };
\node (11) at (5,-3.35) {\large (b) \begin{tabular}{l}  $\widehat{Q}_{i+2}$  needed to compute gradient \\ of  $i+1$ neighborhood. \end{tabular}};
\node (11) at (-3,-7.40) {{\large (c) } \begin{tabular}{l}
$\widehat{Q}_{i+2}$ needs $\widehat{U}_{i+3}$ and overlap count \\ $N_{i+3}$  to form its neighborhood. \end{tabular}};
\node (11) at (5,-7.40) {\large (d) \begin{tabular}{l}
$N_{i+3}$  needs to know  $i+4$ neighborhood,\\ which may need $i+5$ cell volume.
\end{tabular} };
\end{tikzpicture}
\caption{For illustrative purposes, we assume central merging neighborhoods are used for all cells. Each of the neighborhoods is illustrated.}
\label{fig:halocells}
\end{figure}
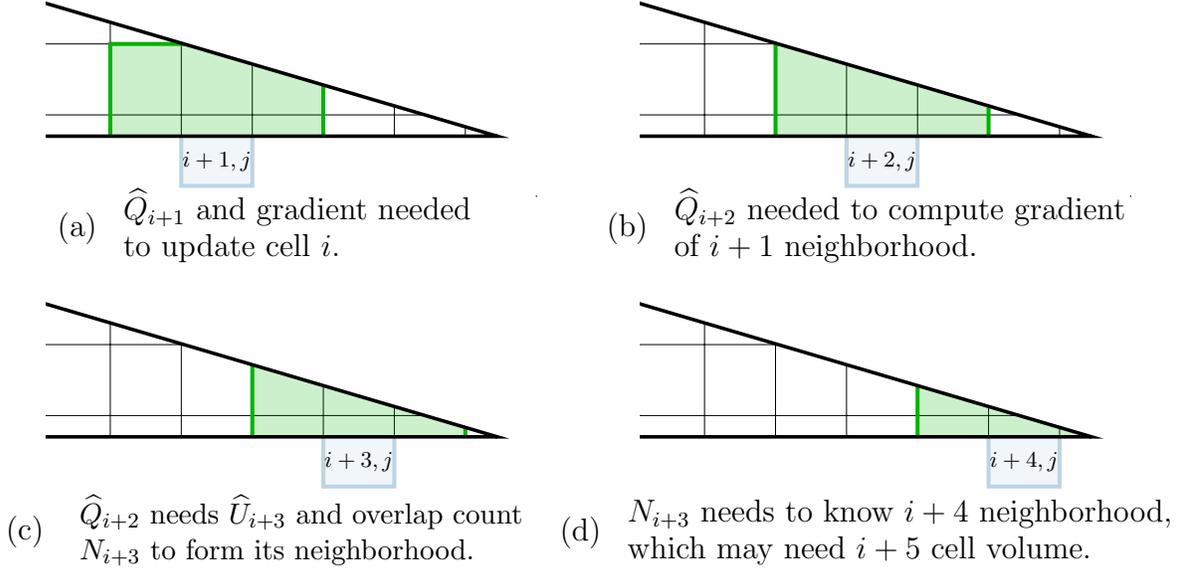

The stencil width for the preprocessing step is determined by how many cells are required in postprocessing, so we count that first using only one  direction.  To update $U_{i,j,k}^{n+1}$ we may need the neighborhood value  ${\widehat{q}_{i+1,j,k}(x_{ijk},y_{ijk},z_{ijk})}$, for example with central merging. (Central merging uses a $3 \times 3 \times 3$ neighborhood around cell $(i,j,k)$). To obtain a slope on the neighborhood associated with cell $(i+1,j,k)$ we need the neighborhood average ${\widehat{Q}_{i+2,j,k}}$. This latter value may have needed ${\widehat{U}_{i+3,j,k}}$, again if using central merging. So the stencil for postprocessing can be up to three cells in every direction. 

To provide  ${\widehat{Q}_{i+2,j,k}}$ for the postprocessing step the term $\widehat{U}_{i+3}/N_{i+3}$ was needed (simplifying subscripts from now on).  $N_{i+3}$ could depend on whether neighborhood $i+4$ overlaps cell $i+3$. So we have to know the details of $i+4$'s neighborhood, which may depends on $i+5$. This is the worst case stencil, depicted in Figure \ref{fig:halocells}. This  determines how many halo cells are needed to preprocess without any additional communication between patches during  single time step. Of course in the simplest case, normal merging will suffice and fewer cells are required. There is also a trade-off between using fewer ghost cells but communicating more often.  However in the example shown later we preprocess using this maximum number of ghost cells without additional communication.

It may happen that there are not enough cells in a stencil for gradient reconstruction on merging neighborhoods. For example, cell $(i+5,j)$ in Figure \ref{fig:halocells} has no $y$ neighbors in its stencil.  For this reason the original finite volume scheme also cannot compute an accurate and well-conditioned $y$ slope. Note that body-fitted meshes at an internal cusp could also have this problem. If a well-conditioned gradient for the base finite volume scheme cannot be reconstructed, then locally the solver and the SRD algorithm drop to first order.

\section{SRD Applications}\label{sec:algext}

The SRD approach was originally applied to hyperbolic conservation laws in two dimensions. 
Here we describe how the SRD methodology is incorporated into a  multi-physics solver  with a more complex time stepping algorithm.  We consider the compressible multi-component Navier-Stokes equations with chemical reactions, and low Mach number models where the flow evolves subject to a constraint. The computational results in section 6 will demonstrate some of these extensions in two geometries of interest.

\subsection{Compressible Reacting Flows}\label{sec:compreact}

For compressible flows, the following conservation equations for mass, species mass fractions, momentum, and energy with a finite rate evaluation of chemistry are solved~\cite{SITARAMAN2021111531}:%
\begin{eqnarray}
&\frac{\partial }{\partial t} \left( \rho \right)  + { \nabla} \cdot \left( \rho \mathbf{u} \right) = 0,\label{eqn:cont}\\
&\frac{\partial }{\partial t} \left( \rho Y_m \right)  + { \nabla} \cdot \left( \rho \mathbf{u} Y_m \right) = -{ \nabla} \cdot \SpeciesFlux_m +  \rho \dot{\omega}_m ,  \qquad m=1,\cdots,N_s \\ %
&\frac{\partial }{\partial t} \left( \rho \mathbf{u} \right)  + { \nabla} \cdot \left( \rho \mathbf{u} \otimes \mathbf{u}  \right) + { \nabla} p = {\nabla} \cdot   \StressTensor, \\ %
& \frac{\partial }{\partial t} \left( \rho E \right)  + { \nabla} \cdot \left( \rho \mathbf{u} E + p \mathbf{u}  \right) = -{ \nabla} \cdot  { \HeatFlux} +
\sum_m h_m(T) \SpeciesFlux_m, \label{eqn:en} \;\; 
\end{eqnarray}
where $\rho$ is the density, $\mathbf{u}$ is the velocity vector, $p$ is the pressure of the mixture, $\rho E$ is the total energy, $N_s$ is the number of species, and $Y_k$ is the mass fraction of the $k$-th species. It is assumed that the species of the gaseous mixture are in thermal equilibrium, at a common temperature $T$. $\StressTensor$ is the viscous stress tensor given by:
\begin{equation}
\StressTensor = \eta \left( \nabla \mathbf{u} + ( \nabla \mathbf{u} )^T\right) + \left( \kappa - \frac{2}{3} \eta\right)  (\nabla \cdot \mathbf{u}) \; \mathbb{I} ,
\end{equation}
where $\eta$ is the shear viscosity and $\kappa$ is the bulk viscosity. The diffusive transport flux of the $m$-th species, $\SpeciesFlux_m$, is approximated using a mixture-averaged diffusion process.  $\dot{\omega}_m$ is the chemical species reaction source term for the $m$-th species. $\HeatFlux$ is the thermal conduction heat flux and $h_m$ is the specific enthalpy of species $m$.

The system  (\ref{eqn:cont})-(\ref{eqn:en}) includes both advective and diffusive fluxes and includes reactions that are potentially stiff on the time scale of advection and diffusion.  The overall integration treats advection and diffusion explicitly and integrates the reactions using a stiff ODE solver.
The time integrator for the compressible equations is a standard second order predictor-corrector approach with an optional fixed point iteration that can be used to tightly couple the reaction and transport terms in the equations, that is similar in spirit to spectral deferred corrections.

Let $I_R$ represent the cell update due to reactions, and $I_{AD}$ represents the advective and diffusive update over a time step. The reactions terms
are computed pointwise by integrating the
ODE representing the reactions. %
Initially at the start of the simulation $I_R$ is obtained by integrating the reaction terms without advection or diffusion.  For later steps, $I_R$ is initialized to the final value from the previous time step; i.e., $I_R^{n,0} = I_R^{n-1,M}$ as defined below.
State redistribution is used to stabilize small cells for advection and diffusion as shown below.

\noindent
Step 1:  Initialize $I_R^{n,0}$ as discussed above.

\noindent
Step 2:  Compute advective and diffusive fluxes, $\mathbf{F}^*_\ell$ at centroids of faces and evaluate preliminary update %
using Eq. (\ref{eq:eb_update}) to obtain
\[
\widehat{U}=
U^n - \frac{1}{V_{i,j,k}}\sum_{
\ell \in \text{faces}}  \mathbf{F}_\ell^*  \cdot \mathbf{n}_\ell A_\ell.
\]
Details of the flux computation are given below.

\noindent Step 3:  Apply state redistribution and compute $I_{AD}^n$ (but omitting the cell indices)
using
\[
I_{AD}^n = \frac{SRD(\widehat{U})- U^n}{\Delta t}.
\]

\noindent Step 4:  Define
\[
U^{n+1,0} = U^n+\Delta t ( I_{AD}^n + I_R^{n,0}).
\]

\noindent
 Steps 5-7 are repeated $M$ times to couple advection, diffusion and reactions.
Formal second order accuracy requires $M$ of at least 1.

\noindent Step 5:
Evaluate fluxes at face centroids using $U^{n+1,m-1}$ and form the flux divergence
\[
I_{AD,i,j,k}^{n+1,m} = - \frac{1}{V_{i,j,k}}\sum_{
\ell \in \text{faces}}  \mathbf{F}_\ell^*  \cdot \mathbf{n}_\ell A_\ell,
\]
and form 
\[
\widehat{U}^m = U^{n+1,m-1}+ \Delta t I_{AD,i,j,k}^{n+1,m} ,
\]
and apply state redistribution to obtain
\[
I_{AD}^{n+1,m} := \frac{SRD(\widehat{U}^m)- U^{n+1,m-1}}{\Delta t}.
\]

\noindent Step 6:
Update the reaction term
\[
I_R^{n,m} = \mathcal{E}_R \left ( U^n,\frac{1}{2} ( I_{AD}^n + I_{AD}^{n+1,m}) \right ),
\]
where $\mathcal{E}_R (U,S)$ evolves the reactions from $t^n$ to $t^{n+1}$ with initial condition $U$ and a constant source term $S$ and then computes the average change due to reactions over the time interval.

\noindent Step 7:
Finally we compute the update
\[
U^{n+1,m}  = U^n + \Delta t\left ( I_R^{n,m}+\frac{1}{2} ( I_{AD}^n + I_{AD}^{n+1,m}) \right ).
\]

The explicit treatment of both advection and diffusion necessitates constraints on the time step for stability.
Using the notation $\mathbf u=(u, v, w)$, we require for each cell in the domain
\[
\frac{(\| \mathbf u \|_{\infty} + c)\Delta t}{h} \le \frac{1}{3}, \;\;
\frac{\eta \Delta t}{\rho h^2 } \le \frac{1}{2d}, \;\;
\]
along with similar restrictions for thermal conductivity and species diffusion.  In the above, $d$ is the spatial dimension of the problem and $h =\mathrm{min}(\Delta x, \Delta y, \Delta z)$.  
We note that in practice, the first constraint is the most restrictive.  Although the diffusive time step restriction scales with $\Delta x^2$, it generally does not constrain the time step, and the convective time step is the most restrictive.%

For the spatial discretization 
we perform a limited extrapolation of the characteristic variables to the centers of the cell faces using a standard van Leer MUSCL limiter.
If one of the cells in the slope computation is completely covered then the one-sided slope in that direction is zeroed.
A two-shock approximation Riemann solver is used at face centroids.  A normal Riemann solver is used to compute the
flux at the EB face.

For the diffusive fluxes, derivatives at cell faces are approximated with a centered second order finite-volume discretization. 
At cut faces with area fraction less than one, the stencil is modified to evaluate the fluxes at face centroids.
At the embedded boundary face, we cast a ray into the interior of the domain and interpolate values onto that ray to compute a one-sided diffusive flux.
By formulating the diffusion in a finite volume form, 
SRD stabilizes this term too (see the remark in section \ref{sec:postprocessing}).

\subsection{Low Mach Number Flow}\label{sec:lowmach}

Next we show how SRD is applied to generic low Mach number flows of the form
\begin{align*}
&\frac{\partial \rho }{\partial t}   + { \nabla} \cdot \left( \rho \mathbf{u} \right) = 0, \\%\label{eqn:lmcont}\\[.05in]
&\frac{\partial \mathbf{u} }{\partial t}   +  \mathbf{u}  \cdot \nabla \mathbf{u}   + \frac{1}{\rho} { \nabla} \pi =\frac{1}{\rho} {\nabla} \cdot   \StressTensor, \\ %
&\frac{\partial \rho \phi}{\partial t} + { \nabla} \cdot \left( \rho \mathbf{u} \phi \right) = R_\phi,\\
\end{align*}%
evolving subject to the constraint
\begin{equation}
    \nabla \cdot \mathbf u = S.
    \label{eq:div_constraint}
\end{equation}
Here $\phi$ corresponds to some set of additional variable such as species mass fractions, enthalpy or an advected scalar, and $R_\phi$ represents additional source terms such as diffusion or reactions. The velocity satisfies an inhomogeneous constraint where $S$ is specified in terms of thermodynamic variables so that the system maintains thermodynamic pressure.  The perturbational pressure, $\pi,$ can be thought of as a Lagrange multiplier that ensures that the evolution of the momentum field is consistent with the constraint. 

For this set of equations, SRD only needs to be applied to the
convective terms, since the reaction terms are evaluated pointwise in each cell, and the diffusive term is treated implicitly.
This results in a time step restriction that requires
\[
\frac{\Delta t~\| \mathbf u\|_{\infty}}{ h} \le \frac{1}{3}
\]
at each point in the domain.
We give an overview of the time advance and indicate at which point SRD is applied.

The equations are discretized using a predictor-corrector method-of-lines 
approach with an approximate projection (see \cite{Musser_MFIX:2021} for details).
In the predictor we define provisional values at time $t^{n+1}$ 
by solving
\begin{eqnarray*}
\rho^{P} &=& \rho^n - \Delta t \; A^n_\rho ,\\
\uvec^{P,*} &=& \uvec^n - \Delta t \; A^n_\uvec + \frac{\Delta t}{\rho} \; D^n - \frac{\Delta t}{\rho} \; G \pi^{\nmh},\\
(\rho  \phi)^{P} &=&  (\rho \phi)^n - \Delta t \; A^n_{\rho \phi} + \Delta t \; R^n_\phi.
\end{eqnarray*}
Here $A_{\rho}^n$ is a discretization of $(\nabla \cdot (\rho \uvec))^n$,
$A_{u}^n$ is a discretization of $(u \cdot \nabla \uvec )^n$, and
$A_{\rho \phi}^n$ is a discretization of $(\nabla \cdot (\rho \phi \uvec))^n$;
$G {\pi}^{n-1/2}$ is an approximation to $\nabla {\pi}^{n-1/2}$, and $D^n$ is an approximation
to $( \nabla \cdot \StressTensor^n)$. After constructing $\uvec^{P,\ast}$, we perform an approximate projection to enforce the constraint~\eqref{eq:div_constraint} for $\uvec^P$ and define  $\pi^{n+1/2,\ast}$. 

In the corrector step, we define the solution at $t^{n+1}$ by solving
\begin{eqnarray*}
\rho^{n+1} &=& \rho^n - \frac{\Delta t}{2} \; (A^n_\rho + A^{P}_\rho) ,\\
\uvec^{n+1,*} &=& \uvec^n - \frac{\Delta t}{2} \; (A^n_\uvec + A^{P}_\uvec)
+ \frac{\Delta t}{2 \rho} \; (D^n + D^{n+1,*})
- \frac{ \Delta t}{\rho} \; G \pi^{\nph,*},\\
(\rho  \phi)^{n+1} &=&  (\rho \phi)^n -  \frac{\Delta t}{2 } \; (A^n_{\rho \phi} + A^{P}_{\rho \phi}) + 
\frac{\Delta t}{2} \; (R^n_{\phi} + R^{n+1}_{\phi}),
\end{eqnarray*}
where the advective terms are defined using the results of the predictor step.
We then apply another approximate nodal projection to enforce the constraint~\eqref{eq:div_constraint} for $\uvec^{n+1}$ and define $\pi^{n+1/2}$.
We emphasize that the SRD step modifies the advective update prior to enforcing the constraint.  The nodal projection, which is applied after the advective and diffusive updates, approximately (to second order in mesh spacing) imposes the constraint regardless of exactly how the advective/diffusive updates were constructed.
Thus, the SRD step does not impact the degree to which the nodal projection enforces the constraint.

In both the predictor and corrector steps outlined above, we construct the advective terms as follows.

\noindent Step 1: Predict, and if necessary limit, the  normal velocity on every cell 
face with non-zero area. For cut faces, this involves all
components of the velocity gradient, since  the cut face centroid is not coordinate-aligned with the cell centroid. 
Given left and right states at each face, solve the Riemann problem with upwinding based on the reconstructed velocities $u_L$ and $u_R$.

\noindent Step 2: Project the normal velocities to satisfy the constraint 
\eqref{eq:div_constraint}.
This defines  $\uvec^{MAC}$.

\noindent Step 3: 
Predict all quantities $s$ on faces, similarly to Step 1. 
Here $s$  represents $\rho,$ $\phi$ and all three velocity components.
At each face solve a Riemann problem with upwinding based on $\uvec^{MAC}.$ 

\noindent Step 4:
For $s = \rho$ or $s = \rho \phi$,  we construct
\begin{eqnarray*}
A_s^*  &=&  \left( (a_{i+\half,j,k} \; \umachi s_{\iph,j,k}- a_{i-\half,j,k} \; \umaclo s_{i-\half,j,k} )  \Delta y \Delta z \right. \\
            &+& (a_{i,j+\half,k} \; \vmachi s_{i,j+\half,k} - a_{i,j-\half,k} \; \vmaclo s_{i,j-\half,k} )  \Delta x \Delta z \\
            &+& \left. (a_{i,j,k+\half} \; \wmachi s_{i,j,k+\half} - a_{i,j,k-\half} \; \wmaclo s_{i,j,k-\half}) \Delta x \Delta y  \right) \; / \; V_{i,j,k}  ,
\end{eqnarray*}
where $a$'s are the area fractions of the cut faces.
If $s$ is one of the velocity components, an additional source, $s_{i,j,k} (D \uvec^{MAC})$, is added to the right hand side.
These formulas do not include an embedded boundary face since the flux is zero there.

\noindent Step 5: Apply state redistribution. To isolate the effect of SRD to the advective terms,  first
define
\[
s^\dag = s^n - \Delta t \; A_s^* \;\;\; \mathrm{for} \;\;\;  s = \rho, \; u, \; \rho \phi.
\]
We then apply SRD to $s^\dag$ to give
\[
s^{SRD} = {SRD}(s^\dag)
\]
and define
\begin{equation*}
A_s = -\frac{s^{SRD}-s^n}{\Delta t}\;\; .
\end{equation*}
Time stepping then continues with the next predictor or corrector step.

As with the other examples, we also apply SRD to the  initial data at the start of the simulation. In this case
the initial data is simply replaced by the redistributed data.

 \section{Computational Examples}\label{sec:compResults}
 
We present 3D results from fluid flow simulations in engineering geometries to demonstrate SRD in  settings with advective, diffusive and reacting updates. 
The examples in sections \ref{sec:pistonbowl} and \ref{sec:cfb} rely respectively on the compressible reacting flow and low Mach number flow models introduced in sections 
\ref{sec:compreact} and \ref{sec:lowmach}. 
 The  SRD algorithm is implemented in AMReX-Hydro as well as a number of other codes including the incompressible flow code {\em IAMR}\footnote{https://github.com/AMReX-Codes/IAMR, 2021}, the
 {\em Pele}\footnote{https://github.com/AMReX-Combustion, 2021} suite of 3 codes, which includes  a low-Mach reacting code, a compressible reacting code, and a physics library, and  {\em MFIX-EXA}\footnote{https://github.com/AMReX-Codes/MFIX-Exa}, a computational fluid dynamics–discrete element model (CFD-DEM) code for low Mach number reacting multiphase flows \cite{Musser_MFIX:2021}.
All of these codes, including the SRD implementation, run on hybrid architectures; the SRD implementation itself runs on both CPUs and GPUs.

\subsection{Compressible Fuel Injection in a Piston-bowl Geometry}\label{sec:pistonbowl}

\begin{figure}
  \centering
  \includegraphics[width=\linewidth]{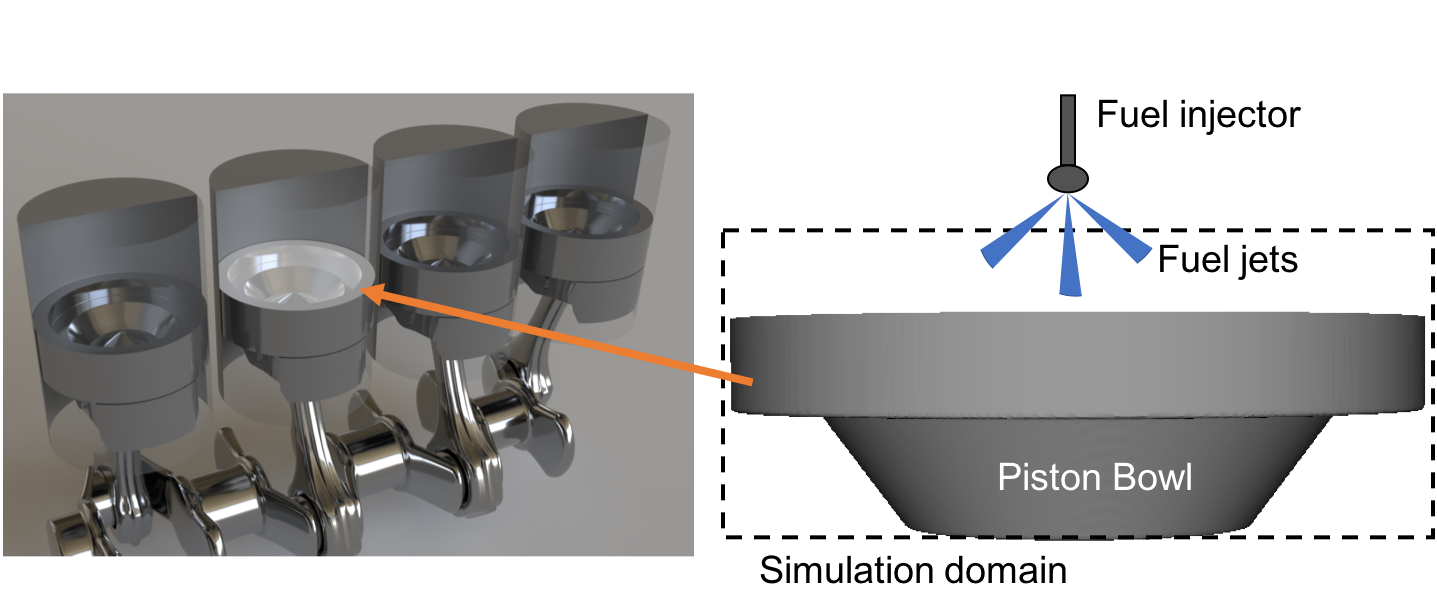}
  \caption{Schematic of the piston bowl geometry, based on a production turbocharged diesel engine.}
  \label{fig:pb-schematic}
\end{figure}

To demonstrate the effectiveness of the state redistribution scheme for compressible flows, {\em PeleC} was used in a simulation of fuel injection in a piston bowl geometry. Complex geometries that concentrate and enhance the flow structure are typical in combustion engines. The geometry used in this work resembles a production turbocharged diesel engine, shown in Figure~\ref{fig:pb-schematic}. The physical domain size is \SI{2.6}{\cm} in the $x$ and $y$ directions and \SI{0.975}{\cm} in the $z$-direction. In compression ignition engines, a multi-hole injector is typically used to inject fuel at the end of the compression stroke. In this work, four discrete gas-phase fuel jets (methane) are injected at the top $z$ boundary into a high temperature oxidizer (oxygen and nitrogen). The inflow velocity conditions for the jets are taken from a turbulent pipe flow precursor simulation. The temperature of the fuel is \SI{300}{\kelvin}. The fuel is injected for \SI{0.00014}{\second} at the start of the simulation. The chemistry kinetics model in this simulation is DRM19 (21 species, 84 reactions) subset of the GRI-Mech 1.2 methane mechanism \cite{Frenklach1995,bell2007numerical}. The ideal gas model equation of state is used to close the equation system.
 
The simulation domain is discretized on the base level with $128 \times 128 \times 48$ cells, leading to a base level grid size of \SI{0.02}{\cm}. The simulation is performed with 2 levels of refinement (effective grid resolution of \SI{0.005}{\cm}) and, when the jets impact the piston bowl side walls, there are 11 million cells on the finest level. This simulation is not possible with the standard flux redistribution scheme as it leads to unphysical solution values (e.g. concentrations less than 0 or $> 1$). The smallest volume fraction allowed in the calculation was $10^{-6}$ to avoid issues with the AMReX EB geometry generation. The SRD algorithm is a negligible portion of the total runtime, typically not exceeding $2\%$.

\begin{figure}
\centering
\begin{subfigure}{.35\textwidth}
  \centering
  \includegraphics[width=\linewidth, trim=0cm 0cm 0cm 0cm, clip=true]{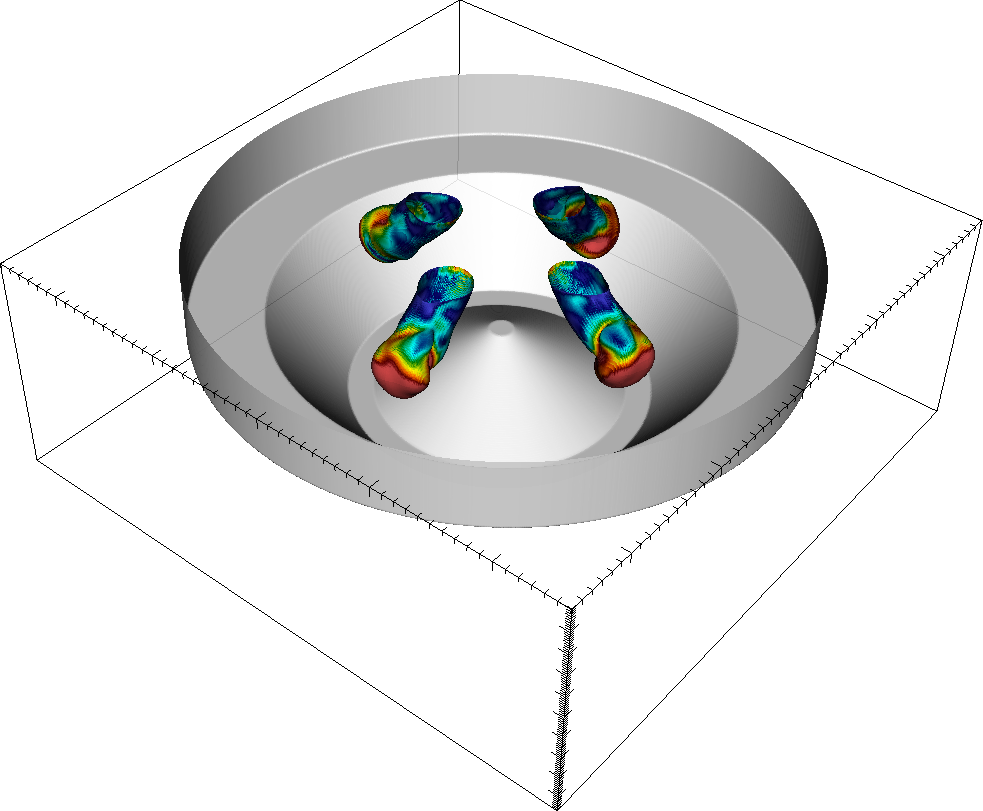}
  \caption{$t=\SI{0.00009}{\second}$.}
  \label{fig:pb-viz-a}
\end{subfigure}\hfill%
\begin{subfigure}{.35\textwidth}
  \centering
  \includegraphics[width=\linewidth, trim=0cm 0cm 0cm 0cm, clip=true]{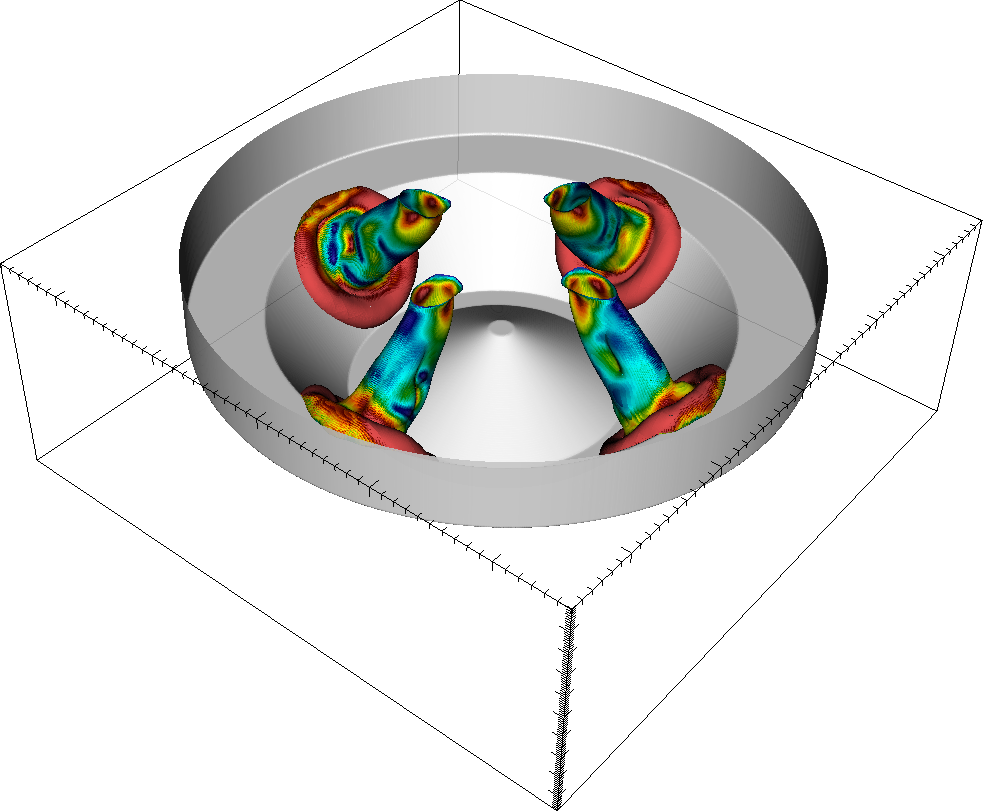}
  \caption{$t=\SI{0.00018}{\second}$.}
  \label{fig:pb-viz-b}
\end{subfigure}\hfill%
\begin{subfigure}{.1\textwidth}
\hspace{1\linewidth}
\end{subfigure}\\%

\begin{subfigure}{.35\textwidth}
  \centering
  \includegraphics[width=\linewidth, trim=0cm 0cm 0cm 0cm, clip=true]{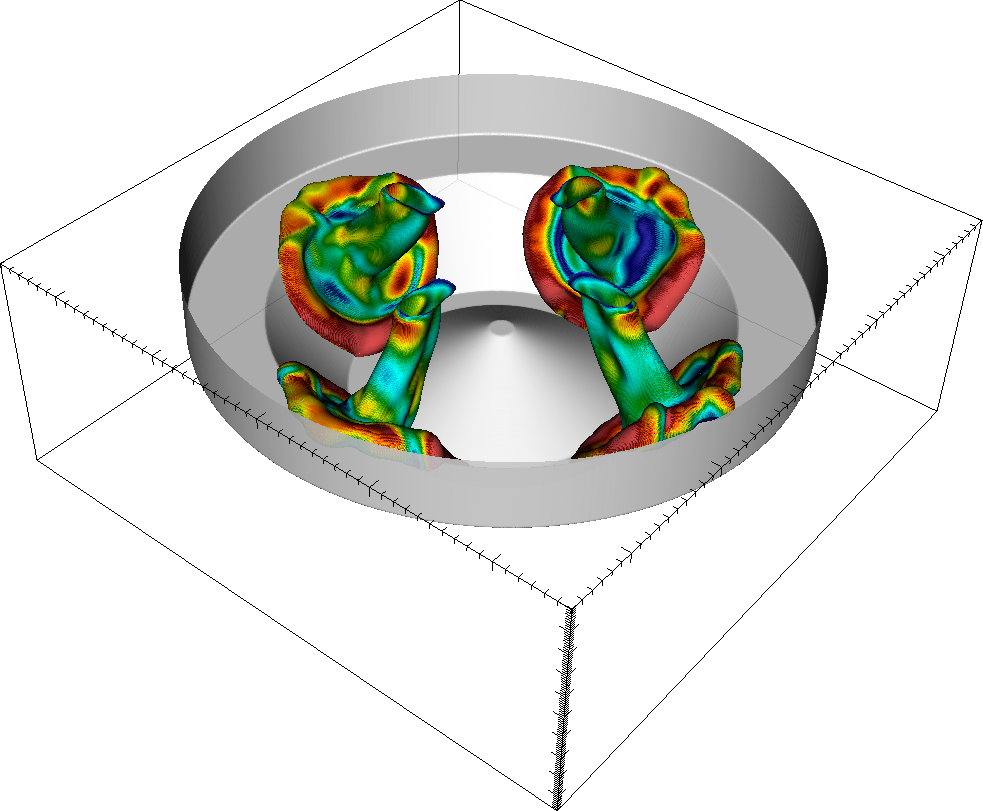}
  \caption{$t=\SI{0.00022}{\second}$.}
  \label{fig:pb-viz-c}
\end{subfigure}\hfill%
\begin{subfigure}{.35\textwidth}
  \centering
  \includegraphics[width=\linewidth, trim=0cm 0cm 0cm 0cm, clip=true]{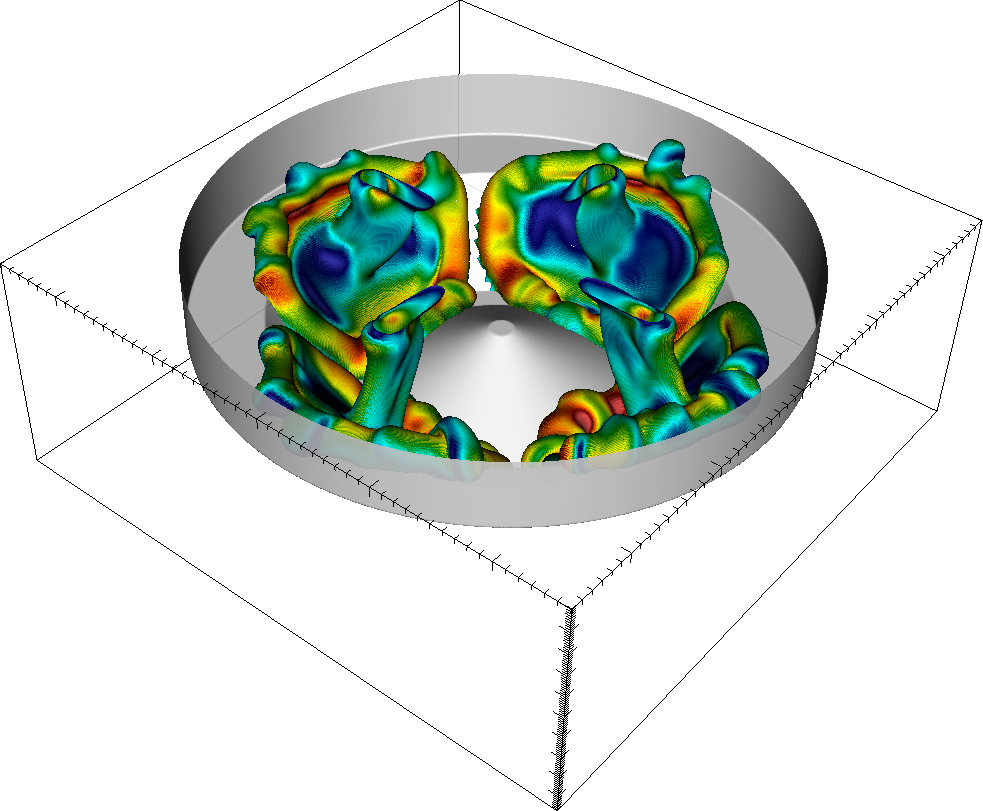}
  \caption{$t=\SI{0.00031}{\second}$.}
  \label{fig:pb-viz-d}
\end{subfigure}\hfill%
\begin{subfigure}{.1\textwidth}
  \centering
  \begin{tikzpicture}
    \node (1) at (0,0) {\includegraphics[width=\linewidth]{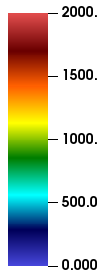}};
    \node (label1) at (0., 2.25) {\small [\SI{}{\cm/\second}]};
  \end{tikzpicture}
  
  \label{fig:}
\end{subfigure}%
\caption{Isosurface of methane mass fraction at 0.1 colored by the magnitude of velocity.}
\label{fig:pb-viz}
\end{figure}

\begin{figure}
\centering
\begin{subfigure}[b]{.49\textwidth}
  \centering
  \includegraphics[width=\linewidth, trim=0cm 3cm 0cm 0cm, clip=true,page=1]{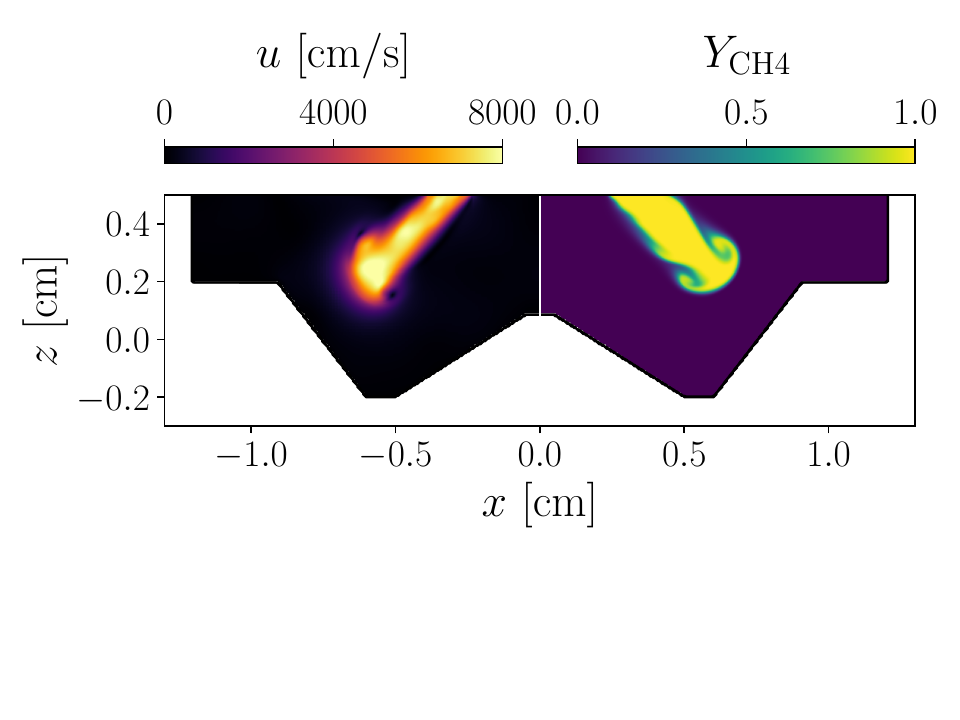}
  \caption{$t=\SI{0.00009}{\second}$.}
  \label{fig:pb-slices-a}
\end{subfigure}\hfill%
\begin{subfigure}[b]{.49\textwidth}
  \centering
  \includegraphics[width=\linewidth, trim=0cm 3cm 0cm 3cm, clip=true,page=1]{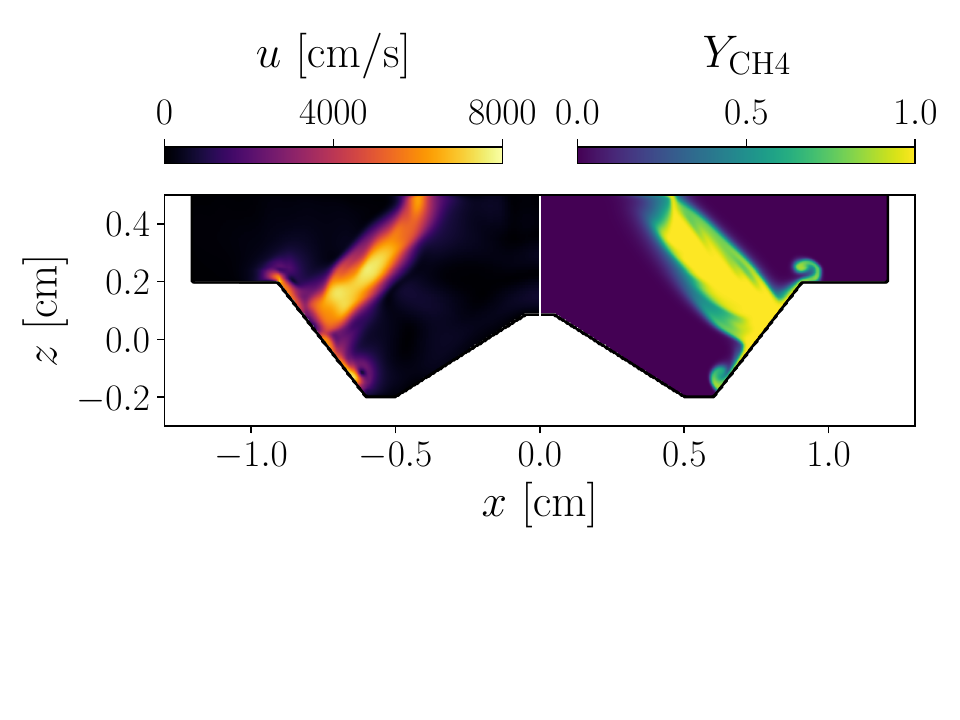}
  \caption{$t=\SI{0.00018}{\second}$.}
  \label{fig:pb-slices-b}
\end{subfigure}\\%
\begin{subfigure}{.49\textwidth}
  \centering
  \includegraphics[width=\linewidth, trim=0cm 3cm 0cm 3cm, clip=true,page=1]{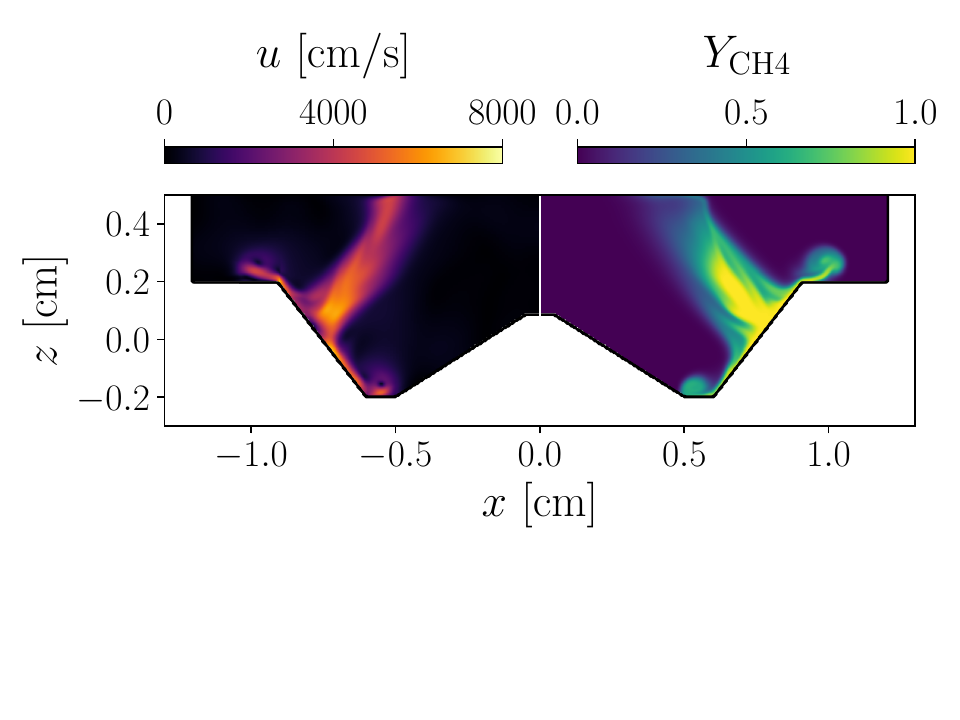}
  \caption{$t=\SI{0.00022}{\second}$.}
  \label{fig:pb-slices-c}
\end{subfigure}\hfill%
\begin{subfigure}{.49\textwidth}
  \centering
  \includegraphics[width=\linewidth, trim=0cm 3cm 0cm 3cm, clip=true,page=1]{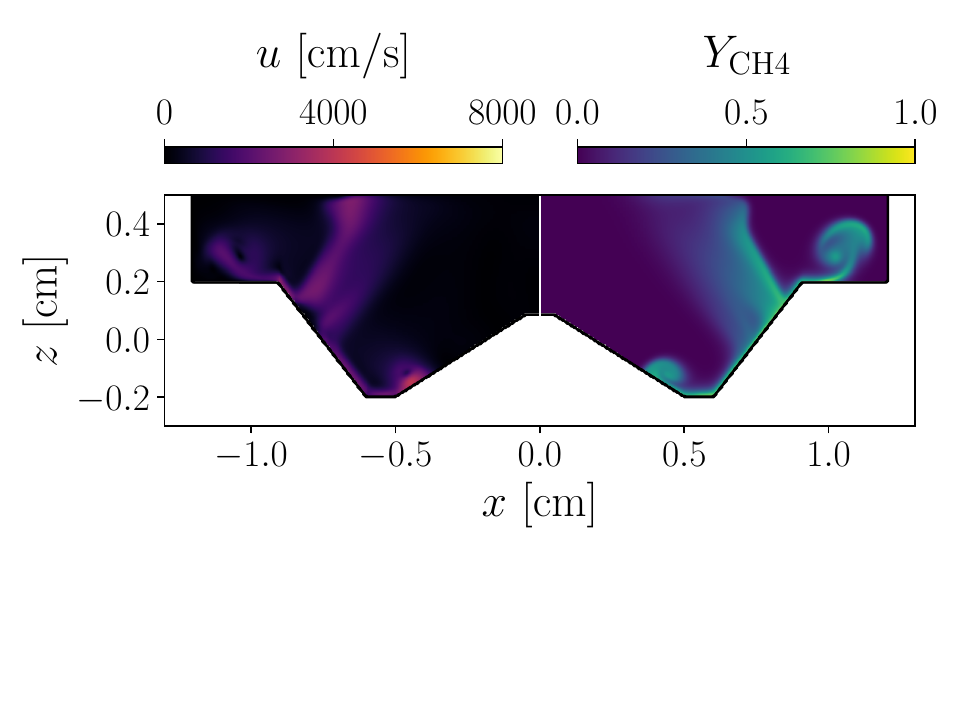}
  \caption{$t=\SI{0.00031}{\second}$.}
  \label{fig:pb-slices-d}
\end{subfigure}%
\caption{Pseudocolors of velocity magnitude (left) and methane mass fraction (right) at $y=\SI{0}{\cm}$.}
\label{fig:pb-slices}
\end{figure}

Simulation results are shown at different times in Figures~\ref{fig:pb-viz} and \ref{fig:pb-slices}. The velocity of the fuel is initially high as the jet penetrates into the domain. The different fluid densities of the oxidizer and fuel lead to characteristic Kelvin-Helmholtz roll-ups at the jet tips, Figures~\ref{fig:pb-viz-a} and \ref{fig:pb-slices-a}. As the fuel jets hit the piston bowl side walls, the fuel is redirected towards the top and the towards the bottom of the bowl, Figures~\ref{fig:pb-viz-b} and \ref{fig:pb-slices-b}. As can be seen in the mass fraction contours in Figure~\ref{fig:pb-slices-b}, the fuel concentration at the inlet has decreased since the fuel stops being injected starting at $t=\SI{0.00014}{\second}$. After impacting the side walls of the piston bowl and interacting with the sides walls, Figures~\ref{fig:pb-viz-c} and \ref{fig:pb-slices-c}, the jets break up and mix into the background fluid, Figures~\ref{fig:pb-viz-d} and \ref{fig:pb-slices-d}.

\subsection{Circulating Fluidized Bed Geometry}\label{sec:cfb}
Fluidization is the phenomenon by which solid particles are converted to a fluid-like phase through the introduction of gas. The resultant mixing of gas and particles provides favorable heat and mass transfer within the system. Such systems are commonly encountered in drying, granulation, coating, heating, and cooling, and over a wide range of industries such as food, agriculture, pharmaceutical, energy and mining. The circulating fluidized bed (CFB) is a type of fluidized bed system that utilizes a recirculating loop for even greater mixing efficiency between the particles and gas. One CFB configuration consists of four main sections – riser, standpipe, loop seal and cyclone. Fluidizing air is introduced primarily from the bottom of the riser. Secondary inflows at lower velocities are also introduced at the bottom of the loop seal and standpipe sections. The gas exits via the cyclone on the top after flowing over the dense particle phase near the bottom. The loop seal returns the particles to the bottom of the riser as they get collected in the standpipe.

\begin{figure}[ht]
\centering
\begin{tikzpicture}
\node (1) at (0,0) {\includegraphics[height=0.4\textheight]{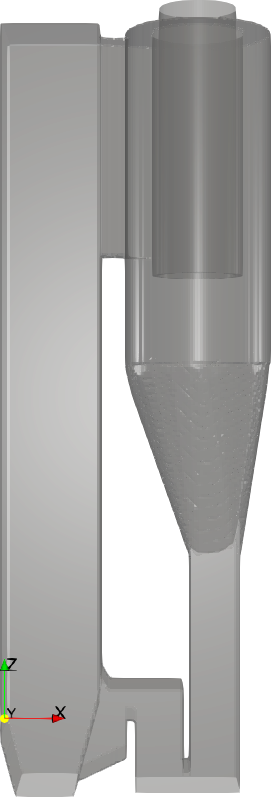}};
\node (riser) at (-2.5,-1) {\small Riser};
\draw[->, line width=1.5] (riser) --++ (1.5,0);
\node (standpipe) at (2.5,-2.3) {\small Standpipe};
\draw[->, line width=1.5] (standpipe) --++ (-1.7,0);
\node (loopseal) at (2.5,-3.1) {\small Loopseal};
\draw[->, line width=1.5] (loopseal) --++ (-2.25,0);
\node (cyclone) at (2.5, 2) {\small Cyclone};
\draw[->, line width=1.5] (cyclone) --++ (-1.75,0);
\node (label1) at (0.75,-4.75) {(a) CFB geometry};

\node (2) at (7,0) {\includegraphics[height=0.4\textheight]{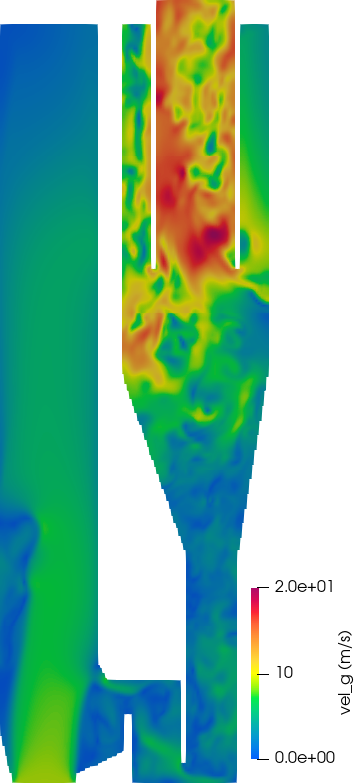}};
\draw [fill=white, color=white] (7.7,-1.8) -- (8.75, -1.8) -- (8.75,-4) -- (7.7, -4) -- (7.7, -1.8);

\node (3) at (9.1,-1.75) {\includegraphics[width=0.04\textwidth, clip, trim= 250 10 80 570]{mfix-figs/square-cfb-velg-wcolorbar_cropped.png}};
\node (0ms) at (9.59, -3.91) {$0$};
\node (10ms) at (9.59, -1.825) {$10$};
\node (20ms) at (9.59, 0.325) {$20$};
\node (ms) at (9.3, 0.85) {\small [\SI{}{\meter/\second}]};
\node (label2) at (7.25,-4.75) {(b) Gas velocity distribution \(y\) = \SI{0.8}{\meter}};

\end{tikzpicture}
\caption{Circulating fluidized bed simulation in MFIX-Exa.}
\label{fig:mfix-cfb}
\end{figure}

MFIX-Exa \cite{Musser_MFIX:2021} was used to simulate such a circulating fluidized bed (CFB) geometry, but without any particles since the focus of this study was only on the fluid phase. (See \cite{Musser_MFIX:2021} for a discussion on how particles are treated in an embedded boundary setting.) In this three-dimensional configuration (Figure~\ref{fig:mfix-cfb}a), the riser, standpipe and loop seal sections have a rectangular cross section, whereas the cyclone has a cylindrical cross section. 

The CFD model for the fluid uses an incompressible flow formulation and employs the SRD algorithm for velocity redistribution. No heat transfer or chemical reactions were considered. The domain of size \SI{3.2}{\meter} in \(x\) direction, \SI{1.6}{\meter} in the \(y\) direction and \SI{8}{\meter} in the \(z\) direction is resolved by a \SI{2.5}{\cm} mesh. Cells whose volume fractions are less than $10^{-6}$ are treated as covered cells in this simulation, due to difficulties in the linear solver that will be the subject of future research. The gas inflow velocity in the riser is increased from \SI{2.5}{\meter/\second} to \SI{9}{\meter/\second} after \SI{1}{\second}. The secondary inflows in the loop seal and standpipes are \SI{0.25}{\meter/\second} and \SI{0.15}{\meter/\second} respectively. Figure~\ref{fig:mfix-cfb}b shows an instantaneous snapshot of the gas velocity magnitude distribution along a vertical slice at \(y\) = \SI{0.8}{\meter} and \(t\) = \SI{1.8}{\second}.

\begin{figure}
    \centering
    \hspace*{-.5in}
    \includegraphics[height=1.8in,trim=0 10 0 0,clip]{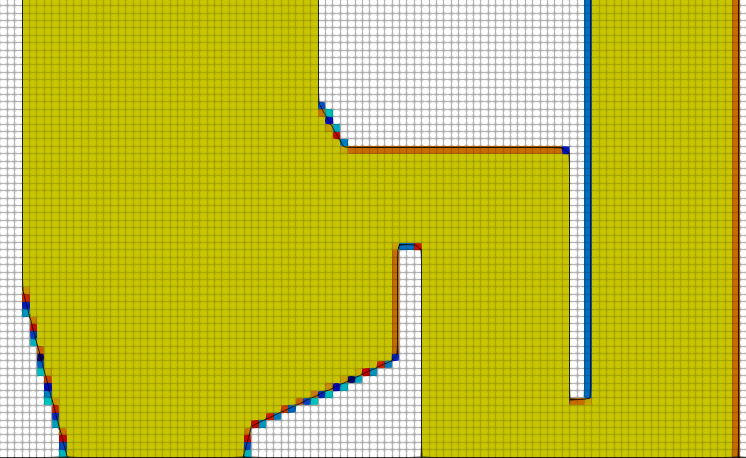}
    \hspace*{.2in}
    \includegraphics[height=1.8in,trim=0 90 0 0,clip]{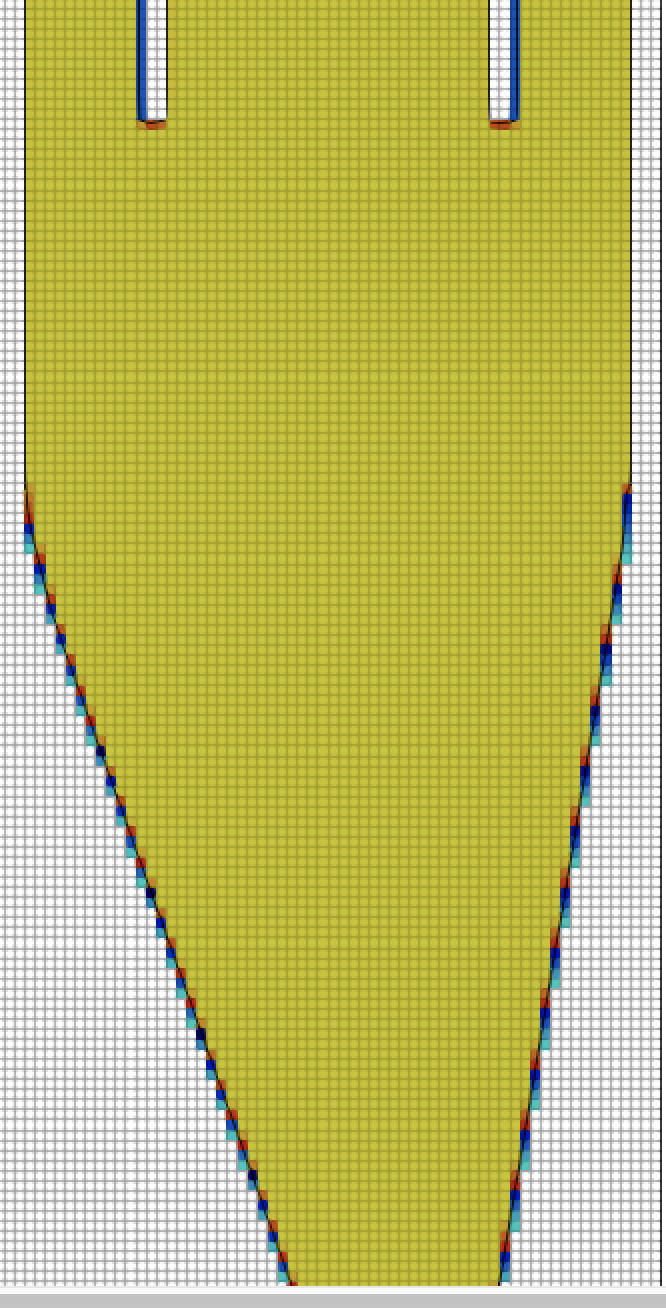}
    \hspace*{.1in}
    \includegraphics[height=1.5in]{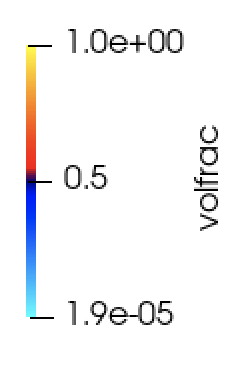}
    \caption{Zoom of grid cross-sections, with colorbar showing cut cell volume fractions.}
    \label{fig:gridFigs}
\end{figure}

The smallest cut cell in this computation had a volume fraction of $9.4 \times 10^{-6}$. Zooms of the grid around the loop seal and the area below the cyclone are shown in Figure \ref{fig:gridFigs}.

 \section{Conclusions and Future Work}\label{sec:conclusions}
We have extended the state redistribution algorithm in several important ways. We introduced a framework to develop new state redistribution methods, and proposed a new variant that activates smoothly with cut cell volume fraction.
The new weighted algorithm is much less sensitive to the merging neighborhood size than the original redistribution algorithm.  Numerical experiments reveal that the weighted algorithm's errors on large neighborhoods is comparable to those on smaller ones.
Next, we showed that SRD can be used in the context of both compressible and incompressible flows, as well as when the PDEs contain reaction and diffusion terms.
All of the presented numerical experiments were obtained using an implementation of SRD in the parallel exascale code AMReX-Hydro.

There are a number of interesting possible extensions of the SRD algorithm presented here.  
For example, one could incorporate knowledge of the local velocity in defining a local weighting scheme for low Mach number flows.  This would be advantageous when the time step is constrained by a velocity $U_{max}$ away from the cut cells.  One could adjust the target volume by the ratio of the local velocity and $U_{max},$ thus potentially eliminating the need to redistribute at all in low-speed regions. For low Mach number flows one could incorporate the pressure gradient as well.
In the context of compressible reacting flows, we would like to formulate a density-weighted SRD algorithm, borrowing ideas used in flux redistribution algorithms \cite{pember:1995}.
We also plan to develop a multi-component slope limiting procedure for the neighborhood slopes to guarantee that chemically reacting species have mass fractions that always sum to 1.

\section*{Acknowledgements}
AG was partially supported by an NSERC (Natural Sciences and Engineering Research Council of Canada) postdoctoral fellowship.  AS and JB were supported by the Exascale Computing Project (17-SC-20-SC), a collaborative effort of the U.S. Department of Energy Office of Science and the National Nuclear Security Administration under contract DE-AC02-05CH11231. MH's work was authored in part by the National Renewable Energy Laboratory, operated by Alliance for Sustainable Energy, LLC, for the U.S. Department of Energy (DOE) under Contract No. DE-AC36-08GO28308. Funding provided by U.S. Department of Energy Office of Science and National Nuclear Security Administration. The views expressed in the article do not necessarily represent the views of the DOE or the U.S. Government. The U.S. Government retains and the publisher, by accepting the article for publication, acknowledges that the U.S. Government retains a nonexclusive, paid-up, irrevocable, worldwide license to publish or reproduce the published form of this work, or allow others to do so, for U.S. Government purposes. This research was supported by the Exascale Computing Project (17-SC-20-SC), a collaborative effort of the U.S. Department of Energy Office of Science and the National Nuclear Security Administration. A portion of the research was performed using computational resources sponsored by the Department of Energy's Office of Energy Efficiency and Renewable Energy and located at the National Renewable Energy Laboratory.

 \bibliography{references}
\bibliographystyle{ieeetr}

\end{document}